\documentclass[letterpaper,11pt]{amsart}
\usepackage[margin=1in]{geometry}

\usepackage{amsmath,amssymb,amsthm,bbm}
\usepackage{dsfont}
\usepackage{mathtools}
\usepackage{enumerate}
\mathtoolsset{showonlyrefs}

\usepackage{hyperref}
\hypersetup{colorlinks=false,pdfborder={0 0 0}}
\usepackage{color}

\usepackage{graphicx}

\newtheorem{theorem}{Theorem}[section]
\newtheorem{proposition}{Proposition}[section]
\newtheorem{definition}{Definition}[section]
\newtheorem{example}{Example}[section]
\newtheorem{remark}{Remark}[section]
\newtheorem{lemma}{Lemma}[section]

\def\.{\hskip.06cm}
\def\ts{\hskip.03cm}

\newcommand\ignore[1]{ }

\def\P{\mathbb{P}}
\def\p{\mathbb{P}}
\DeclareMathOperator{\var}{Var}
\DeclareMathOperator{\e}{\mathbb{E}}

\def\R{\mathbb{R}}
\def\Li{{\rm Li}}

\def\RAc{T_B^c}

\def\supp{{\rm \mbox{supp}}}
\def\L{\mathcal{L}}

\title[Combinatorial QLCLT]{A robust quantitative local central limit theorem with applications to enumerative combinatorics and random combinatorial structures} 
\author{Stephen DeSalvo and Georg Menz}
\date{\today}

\begin{document}
\maketitle
\begin{abstract}
A useful heuristic in the understanding of large random combinatorial structures is the Arratia--Tavare principle, which describes an approximation to the joint distribution of component-sizes using independent random variables. 
The principle outlines conditions under which the total variation distance between the true joint distribution and the approximation should be small, and was successfully exploited by Pittel in the cases of integer partitions and set partitions. 
We provide sufficient conditions for this principle to be true in a general context, valid for certain discrete probability distributions which are \emph{perturbed log-concave}, via a quantitative local central limit theorem. 
We then use it to generalize some classical asymptotic statistics in combinatorial theory, as well as assert some new ones.   
\end{abstract}

{

\noindent \textbf{MSC classes:} 60F05, 60B10, 60C05, 65C50, 05A16 , 05A17, 05A18, 11P82 \\
}

\setcounter{tocdepth}{1}
\tableofcontents

\section{Introduction}

The central limit theorem is a ubiquitous tool which has permeated throughout many diverse fields in mathematics. 
Its universality has found particular importance in enumerative combinatorics, where it paints a broad brush in describing the asymptotic behavior of certain combinatorial sequences, typically using properties of their generating function. 
The Tauberian theorems are of a similar mold, albeit of a more technical nature, with the same goal of describing the growth rates of the coefficients of a generating function. 

More recent in the study of combinatorial sequences, however, is the explicit use of probabilistic tools, which could be considered implicit in many classical approaches. 
One particularly fruitful approach for the so-called ``decomposable" combinatorial structures is expressing the component sizes of a random object, taken uniformly at random from the set of objects with the same \emph{size}, as a joint distribution of globally \emph{dependent} random variables, and approximating this distribution using a joint distribution of \emph{independent} random variables; see~\cite{IPARCS} and the references therein. 
One is then able to approximate various functionals associated with the combinatorial structure via the approximating joint distribution of independent random variables. 

The present paper addresses the phenomena of the local central limit theorem in terms of a smoothness condition, inspired by the applications involving decomposable combinatorial structures. 
We are certainly not the first to attempt such a study, and indeed we are building off the work of our predecessors who have already mapped out large classes of examples for which such local central limit theorems apply. 
Our approach is robust, and is also notably distinct in its assumptions, as we have borrowed a crucial ingredient from enumerative combinatorics, applied in a novel way, meant to expedite the application of our main theorems. 

To quote Hardy and Ramanujan~\cite{HR}:
\begin{quotation}
It is very important, in dealing with such a problem as this, to distinguish clearly the various stages to which we can progress by arguments of a progressively ``deeper" and less elementary character.  The earlier results are naturally (so far as the particular problem is concerned) superseded by the later.  But the more elementary methods are likely to be applicable to other problems in which the more subtle analysis is impractical.
\end{quotation}

To explain our approach, we start with an essential idea from Bender~\cite{Bender}, which is that a local central limit theorem follows from
\[ \mbox{central limit theorem} \quad + \quad \mbox{``smoothness"}. \]
The \emph{central limit theorem} states that for real-valued random variables $X_1, X_2, \ldots$, the partial sums $Y_n = \sum_{i=1}^n X_i$, $n=1,2,\ldots$, are ``close" to a normal distribution in the sense that for all real $x$, we have 
\begin{equation}\label{eq:clt} \sup_x \left| \, \P\left(\frac{Y_n-\e Y_n}{\sqrt{\var(Y_n)}} \leq x\right) - \int_{-\infty}^x \phi(y)\, dy \right| = o(1), \qquad n \to \infty, \end{equation}
where 
\[ \phi(y) = \frac{e^{-\frac{y^2}{2}}}{\sqrt{2\pi}}, \qquad y \in \R, \]
 denotes the density of the standard normal distribution. 
The central limit theorem is robust, in the sense that its statement is the same for continuous or discrete random variables, and is often unchanged when subjected to small perturbations. 
The \emph{local central limit theorem} (LCLT), however, is more delicate; let $\mu = \e Y_n$ and $\sigma^2 = \var(Y_n)$, then for \emph{lattice-valued random variables} only, say with common span $h$, it states that 
\begin{equation}\label{eq:lclt} \sup_{k\in \supp(Y_n)} \left| \frac{\sigma}{h}\ \P(Y_n = k) - \phi\left(\frac{k-\mu}{\sigma}\right) \right| = o(1), \qquad n \to \infty. \end{equation}
When the $o(1)$ in Equation~\eqref{eq:lclt} can be replaced by an expression which is $O(f(n))$, for some function $f(n) \geq 0$ which tends to zero, then we say that $Y_n$ satisfies a \emph{quantitative local central limit theorem} (QLCLT). 

It is often straightforward to verify a central limit theorem, which is why the notion of smoothness will be at the forefront of the analysis. 
Bender~\cite{Bender} performs a detailed analysis of generating functions in order to guarantee smoothness of certain combinatorial sequences of the form $\P(\sum_{i=1}^n X_i = k) = a_{n,k}/\sum_{k=1}^n a_{n,k}$, and shows that a sufficient condition for a local central limit theorem to hold is that the sequence $\{a_{n,k}\}_{n,k \geq 0}$ is \emph{log-concave} in $k$ for each $n$, sometimes called \emph{strongly unimodal}; see Definition~\ref{def:logconcave}. 
Alternatively, if one just assumes the sequence $a_{n,k}$ is unimodal, i.e., that there exists a $k$ such that $a_{n,0} \leq a_{n,1} \leq \ldots \leq a_{n,k-1} \geq a_{n,k} \geq a_{n,k+1} \geq \ldots$, then Bender also obtains a LCLT \emph{except} at a small ball around the central probability, which has the potential to deviate via a sharp peak. 
Such smoothness properties allows one to go from the central limit theorem to the local central limit theorem at the expense of a more detailed analysis of the combinatorial sequence and/or its generating function. 
It is precisely this detailed analysis that we wish to avoid, however, as it is often the case that proving such a condition holds is highly nontrivial for many combinatorial sequences of interest; see for example~\cite{Szekeres2}. 
One advantage this approach has, however, is that it can be applied in contexts where there is dependence between the random variables of interest, whereas we shall assume the random variables are mutually independent. 

A less combinatorial, and more direct probabilistic analysis for proving LCLTs was carried out previously by several authors in~\cite{LCLTMineka, LCLTPetrov, LCLTMukhin, LCLTGamkrelidzeLattice} in the \emph{independent and identically distributed} (i.i.d.) setting. 
We shall be interested solely in the setting when the random variables are discrete, and more specifically lie on some sub-lattice of the natural numbers $\mathbb{Z}$. 
The first trivial observation is that the sub-lattice cannot be, for example, only the even integers, for then taking the sum would only yield an even integer, and our local limit theorem would fail for odd values. 
In the i.i.d.~setting, it was shown in~\cite{LCLTPetrov} that a sufficient condition for a LCLT to hold is for each summand to have a finite absolute third moment and for the range to have a gcd of 1. 
The smoothness in these cases comes from the convolution operation via the sum, which even provides enough regularity for more terms of an asymptotic expansion if enough moments exist, see~\cite[Theorem~2]{LCLTPetrov}. 

One can relax the identically distributed condition, but then we must also modify our sub-lattice condition. 
In fact, it was shown in~\cite{LCLTRozanov} that a necessary condition for a LCLT to hold over integer-valued random variables is 
\begin{equation}\tag{L}\label{lattice:condition} \prod_{k=1}^\infty  \max_{0 \leq m < h} \P(X_k \equiv m\ \mbox{(mod $h$)}) = 0, \qquad \mbox{for any $h \geq 2$}. \end{equation}
There are sufficient conditions given as well, but the resulting theorem is not strong enough to imply a local central limit theorem in the sense of our Theorem~\ref{p_qlclt_for_combinatorics}. 
Other work, in particular in~\cite{LCLTMukhin}, presents alternative sufficient conditions for which a LCLT does hold, which are often difficult to verify or not directly applicable. 
Specifically, a crucial statistic is 
\[ D(X,d) := \inf_{a\in \mathbb{R}} \e \langle (X-a)d\rangle^2, \]
where $X$ is a random variable, $d$ is a real number, and $\langle \alpha \rangle$ denotes the distance of $\alpha$ to its nearest integer. 
Let $T : = \sum_{k=1}^n X_k$ denote the sum of independent random variables $X_1, X_2, \ldots, X_n$, and let $X^\ast  := X - \e X$ denote the symmetrization of a random variable $X$.  The following statistics appear explicitly in~\cite{LCLTMukhin} in the sufficient conditions and rate of convergence for a LCLT to hold: 
\begin{equation}\tag{H}
\label{H:statistics} H(X, d) = \e \langle X^\ast d\rangle^2, \qquad H_n(d) = \sum_{k=1}^n H(X_k, d), \qquad H_n = \inf_{\frac{1}{4} \leq d \leq \frac{1}{2}} H_n(d).  \end{equation}
Aside from the condition that the variance, say $B_n^2$ of the sum tends to infinity, one sufficient condition for a LCLT to hold is that 
\begin{equation}\tag{S}\label{Mukhin:sufficient} B_n^2 = O(H_n), \qquad \mbox{as $n$ tends to $\infty$}.  \end{equation}
There are other sufficient conditions which may be easier to work with, however we would consider verifying any of these conditions to be somewhat arduous, as it involves computing moments of fractional parts of random variables. 

We should strongly note, however, that conditions like~\eqref{Mukhin:sufficient} are typically the most widely applicable, as they can be applied to sums of lattice random variables which have spans of various sizes, subject of course to the necessary condition~\eqref{lattice:condition}. 
Contrast this with the sufficient conditions in~\cite{LCLTGamkrelidzeLattice, Davis, McDonald}, which require all the lattice variables to have a span of 1, as the statistics championed are
\[ q(X) := \sum_{m\in \mathbb{Z}} \min \{\P(X = m), \P(X=m+1)\}, \qquad \delta(X) := \sum_m | \P(X = m) - \P(X=m+1) |, \]
and
\[ s(X) := \sup \{ \P(X = m), \P(X=m+1)\}. \]
The obvious solution of dividing through by the greatest common divisor is not applicable in our setting, since our random variables $X_i$ each have a different span, which very often differs from the greatest common divisor of their sum. 

Another noteworthy approach is by Heinrich~\cite{heinrich1990asymptotic}, where a local central limit theorem is proved for $m$-dependent random fields, along with an asymptotic expansion. 
The approach is the usual one, which is the bulk contribution to the central limit theorem comes from the moments of the sum of the random variables, subject to controlling an oscillatory term. 
The main technical condition is that as $n$ tends to infinity we have 
\[ \frac{\sum_{z \in \mathcal{B}(p)} \e \min_{0 \leq r \leq q-1} \P(Y_z^{(p)} \not\equiv r (mod\ q) \| \mathcal{F}_z^{(p)})}{N_n^2\, \ln(|V_n|^{s-1/2}M_1^{1-s})} \longrightarrow \infty \]
(we do not utilize the notation or definitions in the equation above in the future, and so recounting the full set of definitions would be unnecessarily verbose, hence we refer the interested reader to~\cite{heinrich1990asymptotic}, specifically Section 3). 
This technical condition is of course intimately related to~\eqref{lattice:condition}, and bears a striking resemblance to~\cite[Condition~(II)]{RothSzekeres} for an asymptotic expansion for the number of integer partition numbers subject to certain restrictions. 
One might think of this condition as demanding a concentration on a sublattice which is super-logarithmic in size. 
As in previous instances, verifying that this condition holds is not always straightforward, and in addition it assumes that all summands have a span of 1. 

Stein's method is an accessible approach to obtain central limit theorems for sums of globally dependent random variables, with explicit error bounds~\cite{Stein, chen1975poisson}. 
The dependence is handled via a coupling, and there are many strategies available to demonstrate that the dependence in a particular collection of random variables is tame enough to yield a central limit theorem; see~\cite{chen2010normal}. 
To our knowledge, the majority of theorems which imply a local central limit theorem place assumptions either on the random sum $Y_n$, in the form of a bound on the total variation distance on $(Y_n)^{(i)} = Y_n - X_i$, i.e., 
\[ d_{TV}((Y_n)^{(i)}, (Y_n)^{(i)} + 1); \]
or, more simply, on the individual summands as in
\[ \min\left(\frac{1}{2}, 1 - d_{TV}(X_i, X_i+1)\right). \]
Neither of these approaches is ideal in our setting: either we must know detailed properties of the sum $Y_n$, or to apply the latter condition the span of the random variables must be exactly $1$; see for example~\cite{Rollin, Barbour} and the references therein.  

Let us now describe a class of examples which is particularly resilient to classical assumptions for LCLTs, motivated by the aforementioned applications in enumerative combinatorics. 
Let $X_1, X_2, \ldots$ denote discrete random variables with range in $0, 1, 2, \ldots$, i.e., having a span of 1. 
Let $Y_k := k\, X_k$, $k=1, \ldots, n$; that is, the $Y_k$ are lattice-valued random variables with values in $0, k, 2k, \ldots$, i.e., having a span of $k$. 
Then we would like to show that $\sum_{k=1}^n k\, X_k$ satisfies a \emph{local} central limit theorem under certain conditions which are easy to verify. 
As the random variables $Y_k$ have span at least $k$ for $k \geq 1$, and there does not appear to be any straightforward modification of these conditions that would apply in our setting, we are forced to consider more complicated statistics like those in~\eqref{H:statistics}, which are not easy to verify. 
We thus prefer an extension of the original theory, one that can handle sums of lattice-valued random variables with varying spans which still satisfy enough smoothness to imply a LCLT. 

A partial solution is provided in the recent work of Bogachev~\cite{Bogachev}.  
Building off the extensive history pertaining to limit shapes of combinatorial structures, included is a precise and quantitative local central limit theorem which is broadly applicable to a large class of combinatorial structures. 
As in previous work in the field, the sufficient conditions for a LCLT to hold are placed on the form of the generating function, and in this case one must show that the coefficients of the logarithm of the generating function of the combinatorial sequence satisfy certain summability properties; see also~\cite{odlyzko1995asymptotic}.
A similar and non-overlapping set of conditions also appears in~\cite{Yakubovich}. 
Other authors have provided further generalizations in the context of integer partitions, for example~\cite{Ingham, Goh, NonNegR}, where the set of allowed part sizes is some ``polynomially sparse" set like perfect squares. 
We generalize this treatment to other decomposable combinatorial structures in Section~\ref{section:classes}. 

There are several considerable advantages to a detailed analysis of generating functions, as indeed the same analysis which demonstrates the validity of a LCLT in~\cite{Bogachev} also yields the form of the limit shape. 
Another advantage is the large theoretical exploration concerning such representations, for which there is an extensive literature in combinatorics; see~\cite{Flajolet} and the references therein. 
In particular, multivariate generating functions have recently become more widely accessible to due to the recent work in~\cite{pemantle2013analytic}. 
On the other hand, it is not always a straightforward task to effectively analyze the corresponding generating functions, and for this reason we believe our current approach is advantageous in the fact that it is readily applicable to many problems of interest without extensive auxiliary analysis. 

Our motivation is an approximation principle from the class of decomposable combinatorial structures by Arratia and Tavare~\cite{IPARCS}, where the quality of approximation depends on the rate at which a local central limit theorem (LCLT) holds for a certain sum of \emph{independent} lattice-valued random variables. 
More specifically, we consider random variables $Z_1, Z_2, \ldots$ with distributions given by 
\begin{equation}\label{eq:Z}  \P(Z_i = k) = c_i(x)\, g_i(k)\, x^{i\,k}, \qquad i=1,2,\ldots, \quad k =0,1,2,\ldots,\end{equation}
where $g_i(k)$ is nonnegative,  $x>0$ is such that $\P(Z_i=k)$ is summable, and $c_i(x) = \left(\sum_{k \geq 0} g_i(k) x^{i\,k}\right)^{-1}$ is the normalization factor. 
The joint distribution 
\begin{equation} \label{eq:z} {\bf Z} \equiv {\bf Z}(x) := (Z_1(x), Z_2(x), \ldots, Z_n(x))\end{equation}
 is then an approximation for the true distribution of component sizes 
 \begin{equation}\label{eq:c} {\bf C} \equiv {\bf C}(n) := (C_1(n), \ldots, C_n(n)) \end{equation}
 of a combinatorial structure of size~$n$ determined by the choice $g_i(k)$, where size~$n$ refers to the fact that ${\bf C}$ satisfies $\sum_{i=1}^n i\, C_i(n) = n$. 
For example, with $g_i(k) = 1$, the combinatorial structure is unrestricted integer partitions of size~$n$ and $C_i(n)$ denotes the number of parts of size~$i$; when $g_i(k) = \mathbbm{1}(k\in\{0,1\})$, the combinatorial structure is integer partitions of size~$n$ into distinct parts and $C_i(n) \in \{0,1\}$ similarly denotes the number of parts of size~$i$; when $g_i(k) = \frac{1}{(i!)^k k!}$, we obtain set partitions of size~$n$ and $C_i(n)$ denotes the number of blocks of size~$i$. 

A first observation, one which partially drives this endeavor, is the following: for any $B \subset [n]$, and \emph{for all $x>0$ for which the $Z_i$ are summable},  we have
\begin{equation}\label{sum} \P\left( \sum_{i\in B} i\, Z_i = n \right) = r_B(n)\, x^n \prod_{i\in B} c_i(x), \end{equation}
where $r_B(n)$ corresponds to a combinatorial sequence, the precise meaning of which depends on the choice of $g_i(k)$ and $B$. 
Thus, a LCLT on the left-hand side of Equation~\eqref{sum} provides an independent estimate on the combinatorial sequence $r_B(n)$. 
What makes the sum in the left-hand side of Equation~\eqref{sum} unique is the multiplicative factor $i$, which makes the span of each summand at least $i$, $i=1,2,\ldots,n$, and so as we noted previously many standard approaches either are not applicable, or require extensive calculations in order to certify the validity of the claims. 

In what follows, the notation $\mathcal{L}(X)$ denotes the distribution of the random variable $X$.
\begin{theorem}[\cite{IPARCS}]
Given $g_i(k)$, for each $x>0$ such that the distributions in~\eqref{eq:Z} are summable, and for each $n =1,2,\ldots$, we have 
\[\displaystyle  \mathcal{L}\left((Z_1(x), \ldots, Z_n(x))\middle| \sum_{i=1}^n i\, Z_i = n\right) = \mathcal{L}(C_1(n), \ldots, C_n(n)), \]
where $C_i(n)$ denotes the number of components of size~$i$ in a corresponding random combinatorial structure of size~$n$. 
\end{theorem}

The conditioning event $\{\sum_{i=1}^n i\, Z_i=n\}$ breaks the marginal independence in ${\bf Z}$, and one may ask whether certain subsets of indices $B \subset [n] := \{1,2,\ldots,n\}$ can be chosen such that the distributions $(Z_i)_{i\in B}$ and $(C_i(n))_{i\in B}$ are ``close" in some sense, e.g., Prokhorov distance, or total variation distance. 
Partial answers were provided previously in specific cases, for example Fristedt~\cite{Fristedt} studied the Prokhorov distance for integer partitions, Sachkov~\cite{sachkov1974random} for set partitions. 
In~\cite{IPARCS}, a general principle was outlined for when such approximations should be valid in terms of total variation distance. 
The total variation distance between two \emph{distributions} $\L(X)$ and $\L(Y)$ is given by 
\begin{align}\label{dtv:def} d_{TV}(\L(X), \L(Y)) = \sup_{A \subset \R} \left|\P(X \in A) - \P(Y \in A)\right|, \end{align}
where the $\sup$ is taken over all Borel subsets of $\R$. 
A common abuse of notation which we shall use is to write $d_{TV}(X,Y)$ instead of $d_{TV}(\L(X), \L(Y))$.
Specifically, for any $B \subset [n]$, let 
\begin{equation}\label{eq:TB} T_B := \sum_{i\in B} i\, Z_i, \qquad \mbox{ and } \qquad T \equiv T_n := \sum_{i=1}^n i\, Z_i.  \end{equation}
In~\cite{IPARCS} it is shown that
\begin{align}
\nonumber d_{TV}( (Z_i(x))_{i\in B}, (C_i(n))_{i\in B}) & = d_{TV}( T_B, (T_B | T = n) ) \\
\label{dtv:intro} & = \frac{1}{2} \P(T_B > n) + \frac{1}{2} \sum_{r=0}^n \P(T_B = r) \left|\frac{\P( T_{B^c} = n-r)}{\P(T=n)} - 1 \right|.
 \end{align}
Equation~\eqref{dtv:intro} thus makes the connection between LCLTs and independent process approximations to random combinatorial structures explicit, and also highlights the need for a robust LCLT, since there are many parameters involved in the general setting. 

Equation~\eqref{dtv:intro} was already exploited by Pittel for integer partitions~\cite{PittelShape} and set partitions~\cite{PittelSetPartitions}, largely via an analysis of the corresponding generating function. 
The motivation was to find the \emph{limit shape} of these combinatorial structures.\footnote{A limit shape is a law of large numbers for a graphical representation of a combinatorial structure.} 
\emph{However, the work in~\cite{PittelShape} is more than just a conclusion involving a limit shape.} 
The limit shape is just one statistic of interest for combinatorial structures. 
Aside from the recent treatment in~\cite{Dalal} for concave compositions, recent research appears to be primarily focused on obtaining limit shapes~\cite{NonNegR, Goh, Yakubovich, Bogachev}, which only utilizes a LCLT in the special case of $B = [n]$ (although even the full statement of a LCLT is not necessary, only a sub-exponential lower bound on the probability is required, see for example~\cite{RomikUnpublished}), and the independent process approximation principle appears to have been marginalized. 

We therefore present two main results, one of which is a robust QLCLT under the assumption of \emph{perturbed logconcavity} of the marginal distributions $Z_i(x)$, applicable to sums of the form $\sum_{i\in B} i\, Z_i$ under additional mild assumptions. 
The second is a set of sufficient conditions which provides a resolution to the Arratia-Tavare Principle in a more general context than previously studied. 

In Section~\ref{main:results}, we introduce the necessary definitions and state the main theorems of the paper, which includes sufficient conditions for the Arratia-Tavare principle to hold in a general context.  
In Section~\ref{technical:results}, we state several technical lemmas and use them to prove our main theorems, saving the proofs of the technical lemmas for Section~\ref{section:proofs}. 
In Section~\ref{section:integer:partitions}, we apply our results to integer partitions under various restrictions, demonstrating how to recover many previously obtained theorems and generalizing the work of Pittel to include certain classes of partitions with restrictions.  
In Section~\ref{section:classes}, we apply the theorems to three main classes of combinatorial structures: assemblies, multisets, and selections, and provide motivating examples for each class. 

\section{Main Results}\label{main:results}

\subsection{A general quantitative local central limit theorem} 

\begin{definition}\label{def:logconcave}
A sequence $a_1, a_2, \ldots$ of real numbers is called \emph{log-concave} if $a_n^2 \geq a_{n+1} a_{n-1}$ for all $n \geq 2$.  
We say a discrete random variable $X$ is log-concave if $\P(X=k)$ is a log-concave sequence for $k$ in the support of $X$. 
\end{definition}

\begin{definition}\label{perturbed:unimodal}
We say a random variable $X$ is \emph{perturbed log-concave with constant $C$} if there exists a random variable $Y$ defined on the same support of $X$ such that $\p(Y=k)$ is log-concave in $k$ for each $k$ in the support of $Y$, and there is a positive constant $C_1$ such that
\begin{align}\label{e_def_pertubed_log_concave}
     \frac{1}{C_1} \leq \left|  \frac{\mathbb{P} ( X_k=l )}{\mathbb{P} \left( Y_k =l \right)} \right| \leq C_1;
   \end{align}
\end{definition} 

We state a general quantitative local central limit theorem for the density of the random variable~$\sum_{k=1}^N k\, X_k$, where~$X_k \in \mathbb{N}$ are independent integer-valued random variables.

\begin{theorem}[Quantitative LCLT for decomposable combinatorial structures]\label{p_qlclt_for_combinatorics}
We assume that
\begin{itemize}
\item for~$k \in \mathbb{N}$ the integer-valued random variables $X_{k}\in \mathbb{N}$ are independent;
 \item The random variables~$X_{k}$ are uniformly perturbed log-concave in the sense that there is a constant~$0<C_1< \infty$ such that for any~$k \in \left\{ 1, \ldots, n \right\}$ there is a log-concave random variable~$Y_k$ such that $\supp(Y_k) = \supp(X_k)$ and for all~$l \in \supp(X_k)$ it holds
   \begin{align} 
     \frac{1}{C_1} \leq \left|  \frac{\mathbb{P} ( X_k=l )}{\mathbb{P} \left( Y_k =l \right)} \right| \leq C_1;
   \end{align}
\item the quantities~$\mu$,~$\sigma$,and~$\sigma_{\max}$ are defined by
  \begin{align}\label{e_def_sigma_max}
    \mu := \sum_{k=1}^n k  \mathbb{E} \left[ X_k  \right], \quad  \sigma^2 := \var \left[\sum_{k=1}^n k X_k \right], \quad \mbox{and} \quad \sigma_{\max}^2 := \max_{k} \var \left[ k X_k\right];
  \end{align}
\item for fixed constants~$0< C_2,C_3< \infty $ let~$M$ be defined by
    \begin{align} \label{e_def_set_M}
    M : = \left\{ k \in \mathbb{N} \ \middle| \   \frac{1}{C_2} \leq \frac{\var[k\, X_k]}{\sigma_{\max}^2} \leq C_2 \ \mbox{and} \  \var[X_k] \geq C_3 \right\};
  \end{align}
\end{itemize}
Then there are constants~$0<C_4,C_5, C_6< \infty$ such that 
\begin{align}
\sup_{m}  \left| \frac{\sigma}{h}\, \mathbb{P} \left( \sum_{k=1}^n   k X_k  = m \right) - \phi\left(\frac{m-\mu}{\sigma}\right) \right| \leq  C_6 \frac{\sigma_{\max}}{\sigma}  + C_5\, \sigma \exp \left(- C_4|M| \right),
\end{align}
where the sup is taken over $m \in \supp(\sum_{k=1}^n k\, X_k)$, $h$ is the span of $\sum_{k=1}^n k\, X_k$, 
and $\phi(\cdot)$ is the density of a standard normal random variable. 
\end{theorem}

We emphasize that our local central limit theorem is on the sum of random variables $k X_k$ for $k=1,2,\ldots, n$, where the range of $X_k$ is the set of nonnegative integers. 
For this reason the random variables $k\, X_k$ do not satisfy the classical assumptions traditionally demanded by local central limit theorems. 
When deducing a quanitative central limit theorem, we would also usually expect the assumption that the third moments are bounded (cf.~the Berry-Essen theorem). We have implicitly included this third moment bound by assuming that the random variables are pertubed log-concave. The reason is that for those random variables, one can estimate higher-moments by lower moments (see Lemma~\ref{p_estimating_higher_moments_by_lower_moments}).
We demonstrate the utility of Theorem~\ref{p_qlclt_for_combinatorics} for the purpose of asymptotic enumeration in sections~\ref{section:integer:partitions} and \ref{section:classes}. 

\subsection{The Arratia--Tavare principle}
\label{section:principle}

In~\cite{IPARCS}, a general class of probability distributions is introduced which follows a similar paradigm and connection with combinatorial distributions. 
Let $Z_1, Z_2, \ldots$ denote independent random variables, with distribution given in Equation~\eqref{eq:Z}, viz.,
\begin{equation*} \P(Z_i = k) = c_i(x)\, g_i(k)\, x^{i\,k}, \qquad i=1,2,\ldots, \quad k =0,1,2,\ldots,\end{equation*}
where $g_i(k)$ can be chosen in a multitude of different ways, $x>0$ is any value which makes all of the distributions summable, and $c_i$ is simply the normalizing constant. 

We consider those combinatorial structures which can be described via component sizes. 
For each positive integer $n$, we consider combinatorial objects of size~$n$ with component sizes given in Equation~\eqref{eq:c}, viz., 
\[ {\bf C} \equiv {\bf C}(n) := (C_1(n), \ldots, C_n(n)), \]
where $C_i(n)$ denotes the number of components of size~$i$, for $i=1,2,\ldots, n$; thus, we have $\sum_{i=1}^n i\, C_i(n) = n$. 
By considering the uniform distribution over all combinatorial objects of size~$n$, ${\bf C}(n)$ is a joint distribution of \emph{dependent} random variables. 
Letting $p(n)$ denote the total number of objects of size~$n$, and letting $N(n,{\bf a})$ denote the number of objects of size~$n$ which have component sizes given by ${\bf a} = (a_1, a_2, \ldots, a_n)$, we have 
\[ \P({\bf C}(n) = {\bf a}) = \frac{N(n,{\bf a})}{p(n)}. \]
For example, an integer partition of size~$n$ is a sum of positive integers which sum to $n$.  
In this case, the components are the part sizes, and one may ask, for example, for the distribution of $C_1(n)$ as $n$ tends to infinity. 
By taking $g_i(k) = 1$ in Equation~\eqref{eq:Z}, we have $Z_i$ has a geometric distribution with parameter $1-x^i$, $i=1,2,\ldots$; that $Z_i$ is a good approximation for $C_i(n)$ was also observed in~\cite{Temperley, VershikKerov1977}.  We explore this structure in more detail in Section~\ref{section:integer:partitions}. 

For a second example, a set partition of size~$n$ is a disjoint union of subsets of $\{1, 2, \ldots, n\}$ whose union is $\{1, 2, \ldots, n\}$. 
The components are the subsets, which are more specifically called blocks. 
Taking $g_i(k) = 1/i!^k k!$ in Equation~\eqref{eq:Z}, we have $Z_i$ has a Poisson distribution with parameter $x^i/i!$, $i=1,2,\ldots$; that $Z_i$ serves as an approximation for $C_i(n)$ was observed in~\cite{PittelSetPartitions}. 

Our goal is to maintain as much generality as possible, so that our theorems can be applied as stated to interesting special cases. 

\begin{theorem}[Theorem~3 \cite{IPARCS}]
Suppose for ${\bf a} \in \mathbb{Z}_+^n$ and some functions $f$ and $g_1, g_2, \ldots$ we have 
\[ \P({\bf C}(n) = {\bf a}) = \mathbbm{1}(a_1 + 2a_2 + \ldots + n a_n=n) \frac{f(n)}{p(n)}\prod_{i=1}^n g_i(a_i). \]
Define random variables $Z_1, Z_2, \ldots$ via Equation~\eqref{eq:Z} using the functions $g_1, g_2, \ldots$, i.e.,
\begin{equation*} \P(Z_i = k) = c_i\, g_i(k)\, x^{i\, k}, \qquad k=0, 1,2,\ldots, \end{equation*}
where $c_i = \left(\sum_{k \geq 0} g_i(k) x^{i\,k}\right)^{-1}$ is the normalizing constant. 

For any subset of indices $B \subset \{1, \ldots, n\}$, let $T_B := \sum_{i \in B} i\, Z_i$ denote the weighted sum of random variables over indices $i \in B$.  Let $T \equiv T_n := \sum_{i=1}^n i\, Z_i$.  
Then we have 
\begin{equation}\label{eq:dtv} d_{TV}( (C_i)_{i\in B}, (Z_i)_{i \in B}) = d_{TV}( T_B, (T_B | T_n = n) ). \end{equation}
\end{theorem}
Starting with Equation~\eqref{eq:dtv}, and continuing with the equivalent expression in Equation~\eqref{dtv:intro}, allows us to obtain a bound on the entire process of component-sizes via properties of the unconditioned sums $T_B$, $T_{B^c}$, and $T_n$, namely, 
\begin{equation} d_{TV}(T_B, (T_B | T_n = n)) = \frac{1}{2} \P(T_B > n) + \frac{1}{2} \sum_{r=0}^n \P(T_B = r) \left| \frac{\P(T_{B^c} = n-r)}{\P(T_n = n)} - 1 \right|, \end{equation}
where $T_{B^c} = T_n - T_B$. 

With this more explicit expression for total variation distance, in~\cite{IPARCS} the following heuristic is presented for when the total variation distance ought to be small.  \\

\noindent{\bf Arratia--Tavare principle.}  If
\[ \frac{n-\e T_n}{\sigma_n} \quad \mbox{ is not large} \]
and 
\[ \frac{\e T_B}{\sigma_n} \mbox{ and } \frac{\sigma_B}{\sigma_n}  \mbox{ are small}, \]
then $d_{TV}(C_B, Z_B)$ \emph{should be} small.  \\

Of course, there is still the matter of choosing an appropriate set $B$. 
It was shown in~\cite{PittelShape}, using properties of the generating function, that indeed this principle holds true in the case of integer partitions using $B_- = \{1, \ldots, j_1\}$, with $j_1 = o(\sqrt{n})$, and also for $B_+ = \{j_2, \ldots, n\}$, for $j_2 / \sqrt{n} \to \infty$ and $j_2 \leq (1-\epsilon)\sqrt{n} \log(n) / 2c$ for any $\epsilon>0$.   For set partitions~\cite{PittelSetPartitions}, one can take any $j_1$ and $j_2$ such that $(r-j_1)/\sqrt{\log(n)} \to \infty$ and $(j_2-r)/\sqrt{\log(n)} \to \infty$.  
Thus, the heuristic is very broad since it encompasses an equally broad class of combinatorial structures and joint distributions. 
The approach in both of those cases was to obtain a QLCLT from the respective generating functions, each of which requiring a separate detailed analysis. 
We appeal instead to a property of the random variables $Z_i$, i.e., that of perturbed log-concavity and a variance stabilizing set $M$, which allows us to treat this expression in much greater generality. 

The following definitions are utilized often in what follows. 

\begin{definition}
Define $[n] := \{1,2,\ldots,n\}$.  
Let $B \subset [n]$ denote any subset of $[n]$. 
Let $B^c = [n] \setminus B$, 
\[ T \equiv T_n := \sum_{i=1}^n i\, Z_i, \qquad \qquad T_B := \sum_{i \in B} i\, Z_i, \]
\[ \L(Z_B) := \L((Z_i)_{i\in B}), \qquad \qquad  \L(C_B) := \L((C_i)_{i\in B}), \]
\[ \mu := \sum_{i=1}^n i\, \e Z_i, \qquad  \qquad \mu_B := \sum_{i \in B} i\, \e Z_i, \qquad \qquad \mu_{err} := n-\e T_n,  \]
\[ \sigma^2 := \sum_{k=1}^n \var[k\, Z_k], \qquad \qquad \sigma_B^2 := \sum_{k \in B} \var[k\, Z_k], \qquad \qquad (\sigma_B^c)^2 := \sum_{k \in B^c} \var[k\, Z_k], \]
 \[ \sigma_{\max}^2 := \max_k \var[k\, Z_k], \qquad \qquad \sigma_{B,\max}^2 := \max_{k \in B} \var[k\, Z_k], \qquad \qquad  \sigma_1^2 \equiv \sigma_{1,B}^2 := (\sigma_B^c)^2 / \sigma^2. \] 
Let $\mathcal{N}$ denote a standard normal distribution, and $\mathcal{N}_B$ denotes a normal distribution with mean 0 and variance $\sigma_1^2$.

\end{definition}

We are now ready to introduce sufficient conditions for which the Arratia--Tavare Principle holds, which can be readily verified in many cases of interest. 

\begin{theorem}[Quantitative Arratia--Pittel--Tavare Principle]\label{AT:normal}
Suppose for each $i=1,2,\ldots,$ $g_i(k)$ is a perturbed log-concave sequence in $k$, and define the distributions of $Z_1, Z_2, \ldots$ via Equation~\eqref{eq:Z}, with $x$ chosen such that $\e Z_i^2 < \infty$ for all $i=1,2,\ldots$. 
Fix any constants $0<C_2, C_3 < \infty$, and let 
\begin{align}
M : = \left\{ k \in \mathbb{N} \ | \   \frac{1}{C_2} \leq \left| \frac{\var[k Z_k]}{ \sigma_{\max}^2} \right| \leq C_2 \mbox{ and } \var[Z_k] \geq C_3 \right\};
\end{align}
and for any set $B \subset [n]$, define $M_B := M \cap B$, and  $M_B^c := M \cap B^c$. 
Finally, define
\[ \L(C_B) \equiv \L(C_B(m)) := \L(Z_B | T_n=m), \qquad  \mbox{for $m \in \supp(T_n)$}. \]
Then 
  \begin{enumerate}
\item If $T_B$ satisfies a \emph{central limit theorem} with rate $O(r_B)$; or, more strongly, if $\sigma_{B,\max}/\sigma_B \to 0$ and $\sigma_B e^{-C\, |M_B|} \to 0$, then for any $m \in \supp(T_n)$, we have 
\begin{align}\label{dtv:B} 
d_{TV}(C_B, Z_B) & = d_{TV}(\mathcal{N}_B, \mathcal{N}) + O\biggl(\min\left(r_B,\frac{\sigma_{B,\max}}{\sigma_B} + \sigma_B e^{-C\, |M_B|}\right) + \exp\left(- \frac{n-\mu_B}{\sigma_B} \right) +  \frac{\sigma_B\ \mu_{err}}{ (\sigma_B^c)^2} + \frac{\mu_{err}^2}{(\sigma_B^c)^2} + \\
& \qquad \frac{\sigma_{\max}}{\sigma} + \sigma\, e^{-C_3 |M|} + \frac{\sigma_{B^c,\max}}{\sigma_{B^c}} + \sigma_{B^c}\, e^{-C_3 |M_B^c|}\biggr).
\end{align}
\item If $T_B$ is such that $\sigma_1 \to 1$, then we have 
\begin{align}\label{dtv:B:zero} d_{TV}(C_B, Z_B) & = O\biggl( \left(1- \sigma_1^2\right) +\min\left(r_B,\frac{\sigma_{B,\max}}{\sigma_B} + \sigma_B e^{-C\, |M_B|}\right)+\exp\left(- \frac{n-\mu_B}{\sigma_B} \right) +  \frac{\sigma_B\ \mu_{err}}{ (\sigma_B^c)^2} + \frac{\mu_{err}^2}{(\sigma_B^c)^2} + \\
& \qquad \frac{\sigma_{\max}}{\sigma} + \sigma\, e^{-C_3 |M|} + \frac{\sigma_{B^c,\max}}{\sigma_{B^c}} + \sigma_{B^c}\, e^{-C_3 |M_B^c|}\biggr).
\end{align}
\end{enumerate}
\end{theorem}

We defer the proof until Section~\ref{AT:proof}, after we have established the relevant quantitative local central limit theorem. 
Note that Theorem~\ref{AT:normal} is robust in the sense that some of the $Z_i$ are allowed to be point masses at $0$, which is why we also need the assumption that $m \in \supp(T_n)$. 
We give examples in sections~\ref{section:integer:partitions}  and~\ref{section:classes}.  
The following is a simpler, non-quantitative version of Theorem~\ref{AT:normal}. 

\begin{theorem}[Sufficient conditions for the Arratia--Pittel--Tavare principle]\label{AT:principle} 
Recall the definitions in Theorem~\ref{AT:normal}. 
Suppose for some $0<C_3<\infty$ we have
\[ \frac{\mu_{err}}{\sigma_{B^c}},\ \ \frac{\sigma_B}{n-\mu_B}, \ \ \frac{\sigma_{\max}}{\sigma}, \ \ \sigma\, e^{-C_3 |M|}, \ \ \frac{\sigma_{B^c, \max}}{\sigma_{B^c}}, \ \ \sigma_{B^c}\, e^{-C_3 |M_{B^c}|} \to 0; \] \vskip .025in 
\[ \sigma_B = O(\sigma_{B^c});  \qquad \qquad  \sigma_1^2 \to 1; \] 
then $d_{TV}( C_B, Z_B ) \to 0$.  
\end{theorem}

The case when $\displaystyle \lim \sigma_1  \in (0,1)$ is interesting in its own right, i.e., when the total variation distance tends to a nontrivial explicit limit; see Example~\ref{dtv:root:two} for an application which exploits this case.  
When $\sigma_1 \to 0$ the argument fails altogether, which corresponds to a greedy strategy of including too many components, i.e., taking the set $B$ too large. 

\section{Technical results}
\label{technical:results}

\subsection{Technical lemmas}

Fundamental for our proofs is the observation that as soon as a random variable~$X$ is normalized one has a uniform control of the corresponding characteristic function~$\phi(\xi) = \mathbb{E} \left[ \exp (i X \xi) \right]$ close to~$\xi=0$. 
This control is obtained via an assumption on the absolute third moment, which we quantify in Lemma~\ref{p_taylor_characteristic_function} below.

\begin{lemma}[Control of normalized characteristic functions close to~$0$]\label{p_taylor_characteristic_function}
Assume that the random variable~$X$ is normalized in the sense that
\begin{align*}
  \mathbb{E}[X] =0 \qquad \mbox{and} \qquad \var [X] =1,
\end{align*}
and the absolute third moment is bounded i.e.
\begin{align*}
  \e |X|^3 < \infty.
\end{align*}
Then for each $\delta>0$, there is a finite, positive constant $C \equiv C(\delta)$, such the function~$h(\xi)= - \log \mathbb{E}\left[ \exp \left(iX \xi \right)\right]$ satisfies for all~$|\xi| \leq \delta$ the estimate
\begin{align}\label{e_control_characteristic_taylor}
 \left| h (\xi) - \frac{1}{2} |\xi|^2 \right| \leq C \e |X|^3\,  |\xi^3|.
\end{align}
\end{lemma}

Equally essential is a control on the characteristic function away from 0 (see Lemma~\ref{p_control_characteristic_fnction_intermediate_and_large_values} from below). 
This is where the majority of the technical work lies, as well as the novelty in our approach, since we must assume some type of structure which eliminates a type of periodicity. 
In our case, we have chosen to place a rather mild and easily satisfied technical condition on the probability mass functions of the random variables. 
In addition, the argument for deducing Lemma~\ref{p_control_characteristic_fnction_intermediate_and_large_values} is quite robust, as we also get uniform control on the characteristic function in other situations which are typically specialized to particular applications involving a detailed analysis of the generating function. 

The perturbed log-concave assumption is critical in the proof of the following lemma, Lemma~\ref{p_control_characteristic_fnction_intermediate_and_large_values} below, which gives enough control to prevent gaps in the local central limit theorem. 
Its proof is inspired by the proof in~\cite[Lemma~3.4]{Menz_Otto} for real-valued random variables. 

\begin{lemma}[Control of characteristic functions away from 0] \label{p_control_characteristic_fnction_intermediate_and_large_values}
We consider a random variable~$X$ on a countable probability space. We assume that the values of random variable~$X$ are a.s.~on a one-dimensional lattice with distance smaller than~$1$ between the lattice points. More precisely, we assume that there exists bounded constants~$0 \leq c_1, c_2 \leq 1$ such that
\begin{align}\label{e_lattice_random_variable}
  \mathbb{P} \left( X= c_1 k + c_2, \ k \in \mathbb{Z}  \right) =1.
\end{align}
We assume that the first moment of~$X$ is bounded from above, i.e.
\begin{align*}
  \mathbb{E} |X| \leq C_2< \infty,
\end{align*}
and that the variance is bounded from below, i.e. there is~$C_3 < \infty$ such that
\begin{align*}
  \var (X) \geq \frac{1}{C_3}.
\end{align*}
Additionally, we assume that there is a constant~$0<C_4 < \infty$ such that the random variable~$X$ is perturbed log-concave with constant $C_4$.
Let~$h(\xi)= \mathbb{E}\left[\exp (i X \xi) \right]$ denote the characteristic function of~$X$. Then for any~$\delta>0$ there exists a constant~$\lambda<1$, depending only on~$\delta$,~$C_1$,~$C_2$,~$C_3$ and on~$C_4$ such that for all~$\delta\leq |\xi| \leq 2 \pi - \delta$
\begin{align}\label{e_uniform_bound_characteristic_function_away_from_0}
  |h(\xi)| \leq \lambda.
\end{align}
\end{lemma}

Finally, the fact that our random variables are assumed to be perturbed log-concave 
allows us to estimate higher central moments by lower central moments. 
The statement of Lemma~\ref{p_estimating_higher_moments_by_lower_moments} below is well-known in the context of real-valued random variables (see, e.g.,~\cite[Section 3]{BoMA11} and the references therein.)

\begin{lemma}[Reverse H\"older inequality for central moments] \label{p_estimating_higher_moments_by_lower_moments}
  Assume that the random variable~$X \in \mathbb{N}$ is perturbed log-concave (see Definition~\ref{perturbed:unimodal}). Then there is a universal constant~$0< C < \infty$ such that
  \begin{align*}
    \mathbb{E} \left[ \left| X - \mathbb{E} [X] \right|^3 \right] \leq C C_1^{1+\frac{3}{2}} \left( \var(X) \right)^{\frac{3}{2}},
  \end{align*}
where the constant~$C_1$ is given by~\eqref{e_def_pertubed_log_concave}. More generally, it holds that there is a universal constant~$0< C < \infty$ such that for any~$ p,q \in \mathbb{N}$satisfying~$1 \leq q \leq p $ it holds
\begin{align}\label{e_inverse_hoelder}
    \left( \mathbb{E} \left[ \left| X - \mathbb{E} [X] \right|^p \right] \right)^{\frac{1}{p}} \leq C C_1^{\frac{1}{p}+ \frac{1}{q}} p  \left( \mathbb{E} \left[ \left| X - \mathbb{E} [X] \right|^q \right] \right)^{\frac{1}{q}}.
\end{align}
\end{lemma}

These lemmas are proved in Section~\ref{section:proofs}. 
We now have all of the ingredients to prove our main theorem. 

\begin{proof}[Proof of Theorem~\ref{p_qlclt_for_combinatorics}]
Let $W : = \sum_{k=1}^n k X_k $, and for now assume that $h = 1$, i.e., the span of $W$ is 1.
Using the inversion theorem, we have 
  \begin{align}
    \mathbb{P} \left( W = m \right) & = \frac{1}{2 \pi} \int_{- \pi}^{\pi}
  \mathbb{E}\left[ \exp \left(i \left( W - m \right) \xi \right) \right] d \xi \\
    & =   \frac{1}{2 \pi} \int_{- \pi}^{\pi}   \mathbb{E}\left[ \exp \left(i  \sum_{k=1}^n \left( X_k - \mathbb{E}\left[X_k \right] \right) \xi \right) \right] d \xi \\
    &  = \frac{1}{2 \pi} \int_{- \pi}^{\pi}   \exp\left(i\, \xi (\mu-m)\right) \mathbb{E}\left[ \exp \left(i  \sum_{k=1}^n  \frac{1}{\sigma_k}\left( kX_k - \mathbb{E}\left[kX_k \right] \right) \sigma_k \xi \right) \right] d \xi \\
    & =  \frac{1}{2 \pi} \int_{- \pi}^{\pi}   \exp\left(i\, \xi (\mu-m)\right) \mathbb{E}\left[ \exp \left(i  \sum_{k=1}^n \hat X_k  \sigma_k  \xi \right) \right] d \xi , \label{e_rewriting_P}
  \end{align}
where we used the notation
\begin{align}
\mu = \sum_{i=1}^n \e{iX_i}, \qquad \mbox{ and } \qquad  \sigma_k^2 = \var \left(k X_k \right) \qquad \mbox{ and } \qquad   \hat X_k =  \frac{1}{\sigma_k}\left( kX_k - \mathbb{E}\left[kX_k \right]\right) . \label{e_normalization}
\end{align}
Next, let $\phi_{\mu,\sigma}(t) = e^{-(t-\mu)^2/2\sigma^2}/\sqrt{2\pi \sigma^2}$ denote the density function of a normal random variable. 
Again by the inversion theorem, we have 
\begin{align}
\label{normal:equation} \phi_{\mu, \sigma}(m) & = \frac{1}{2\pi} \int_{-\pi}^{\pi} e^{-i \xi (m-\mu)} e^{ - \xi^2 \sigma^2 / 2} d\xi. 
\end{align}
We now subtract~\eqref{normal:equation} and~\eqref{e_rewriting_P} and obtain uniformly in $m \in \supp(W)$ 
\begin{equation}\label{uniform:bound} \left| \mathbb{P} \left(W = m \right)  - \phi_{\mu,\sigma}(m)\right| \leq 
\frac{1}{2\pi} \int_{-\pi}^{\pi}  \left| \mathbb{E}\left[ \exp \left(i  \sum_{k=1}^n \hat X_k  \sigma_k  \xi \right) \right] -  \exp \left( - \frac{\sigma^2 |\xi|^2}{2} \right)\right| d \xi .
\end{equation}

We choose~$\delta>0$ small but fixed and we want to recall the definition~\eqref{e_def_sigma_max} of~$\sigma_{\max}$. We split up the integral on the right hand side of~\eqref{uniform:bound} into two parts; i.e.
\begin{align}
& \frac{1}{2\pi} \int_{-\pi}^{\pi}  \left| \mathbb{E}\left[ \exp \left(i  \sum_{k=1}^n \hat X_k  \sigma_k  \xi \right) \right] -  \exp \left( - \frac{\sigma^2 |\xi|^2}{2} \right)\right| d \xi \\
 & \qquad + \frac{1}{2\pi}  \int_{\frac{\delta}{\sigma_{\max}} \leq |\xi| }  \left| \mathbb{E}\left[ \exp \left(i  \sum_{k=1}^n \hat X_k  \sigma_k  \xi \right) \right] - \exp \left( - \frac{\sigma^2 |\xi|^2}{2} \right) \right| d \xi  \\
& \quad =: T_1 + T_2. \\
\end{align}

The term~$T_1$ on the right hand side of the last estimate is called the inner integral. The term~$T_2$ is called the outer integral. We will estimate each term separately. We start with estimating the term~$T_1$. We rewrite the first integrand of the term~$T_1$ according to the definition of $h$ in Lemma~\ref{p_taylor_characteristic_function}, viz., 
\begin{align*}
 \mathbb{E}\left[ \exp \left(i  \sum_{k=1}^n \hat X_k  \sigma_k  \xi \right) \right]
& =  \mathbb{E}\left[ \exp \left( \sum_{k=1}^n h \left( \sigma_k  \xi \right)\right) \right] .
\end{align*}
We want to note that by Lemma~\ref{p_estimating_higher_moments_by_lower_moments} we have that there is a universal constant~$0< C< \infty$ such that
\begin{align*}
  \varrho_k &:= \mathbb{E} \left[ | kX_k - \mathbb{E}[kX_k^3]|^3 \right] \\
 & \leq C k^3 \left( \var (X_k) \right)^{\frac{3}{2}} \\
 & = C  \left( \var (kX_k) \right)^{\frac{3}{2}} \\
  & = C \sigma_k^{\frac{3}{2}}.
\end{align*}
Hence, by applying Lemma~\ref{p_taylor_characteristic_function} we get
\begin{align}
  \left| \sum_{k=1}^n h_k \left( \sigma_k \xi \right)  - \frac{1}{2} \sigma^2 \xi^2\right| & \leq \left| \sum_{k=1}^n h_k \left( \sigma_k \xi \right)  - \frac{1}{2} \sum_{k=1}^n \sigma_k^2 \xi^2 \right| 
 \leq   \sum_{k=1}^n  \left|h_k \left( \sigma_k \xi \right) - \frac{1}{2} \sigma_k^2 \xi^2 \right| \\
& \leq  C \sum_{k=1}^n \rho_k  \left| \xi \right|^3  
 \leq  C \left(\frac{\sum_{k=1}^n \rho_k }{\sigma^2}\right)  \sigma^2 \left| \xi \right|^3  \\
  & \leq  C  \sigma_{\max}|\xi| \   \sigma^2 \left| \xi \right|^2 
   \leq  C \delta  \sigma^2 \left| \xi \right|^2 , \label{e_first_bound_by_taylor_1}
\end{align}
where we used that in the inner integral~$T_1$ it holds~$ \sigma_{\max} |\xi| \leq \delta$. By choosing~$\delta$ small enough, it follows from~\eqref{e_first_bound_by_taylor_1} that
\begin{align*}
  \mbox{Re} \sum_{k=1}^n h_k \left( \sigma_k \xi \right)  \geq  \frac{1}{4} \sigma^2 \left| \xi \right|^2.
\end{align*}
Using now the Lipschitz continuity of~$\mathbb{C} \ni y \mapsto \exp (y)$ on~$\mbox{Re} (y) \leq - \frac{1}{4} \sigma^2 \left| \xi \right|^2$ with constant~$\exp(- \frac{1}{4} \sigma^2 \left| \xi \right|^2)$ yields the estimate
\begin{align*}
  \left| \exp \left( \sum_{k=1}^n h_k \left( \sigma_k \xi \right) \right) - \exp\left(  -  \sigma^2 \xi^2 \right) \right| & \leq C \sum_{k=1}^n \rho_k  |\xi|^3  \exp\left(- \frac{1}{2}  \sigma^2 \xi^2 \right). 
\end{align*}
Plugging this estimate into the term~$T_1$ yields
\begin{align*}
T_1 & = \frac{1}{2 \pi} \int_{|\xi|\leq \frac{\delta}{\sigma_{\max}} }  \left| \exp \left(   \sum_{k=1}^n h_k\left( \sigma_k \xi \right) \right) -  \exp \left( - \frac{\sigma^2 |\xi|^2}{2} \right) \right| d \xi  \\
& \leq  C \left(\sum_{k=1}^n \rho_k \right) \int_{|\xi| \leq \frac{\delta}{\sigma_{\max}}} |\xi|^3   \exp\left(- \frac{1}{4}  \sigma^2 \xi^2 \right) d\xi\\
& =  C\, \frac{\left(\sum_{k=1}^n \rho_k\right)}{\sigma^4} \int_{|\xi| \leq \frac{\delta}{\sigma_{\max}}} |\sigma\, \xi|^3   \exp\left(- \frac{1}{4}  \sigma^2 \xi^2 \right) d(\sigma\, \xi)\\
& = C \frac{\left(\sum_{k=1}^n \rho_k\right)}{\sigma^4}   \int_{|\xi| \leq \frac{\delta}{\sigma_{\max}} \sigma} |\xi|^3   \exp\left(- \frac{1}{4}  \xi^2 \right) d\xi \\
& \leq C \frac{\left(\sum_{k=1}^n \rho_k\right)}{\sigma^4}  \\
& \leq C \frac{\sigma_{\max}}{\sigma^2}.
\end{align*}

Now, let us estimate the outer integral given by the term~$T_2$. We start with applying the triangle inequality and get
\begin{align}
  T_2 & \leq \frac{1}{2\pi} \int_{\frac{\delta}{\sigma_{\max}} \leq |\xi|  } \left|   \mathbb{E}\left[ \exp \left(i  \sum_{k=1}^n \hat X_k  \sigma_k \xi \right) \right] d \xi -  \int_{ \frac{\delta}{\sigma_{\max}} \leq |\xi| } \exp \left( - \frac{\sigma^2 |\xi|^2}{2} \right) \right| d \xi  \\
  & \leq \frac{1}{2\pi}  \int_{\frac{\delta}{\sigma_{\max}} \leq |\xi|  }  \left| \mathbb{E}\left[ \exp \left(i  \sum_{k=1}^n \hat X_k  \sigma_k \xi \right) \right]\right|  d \xi  +  \frac{1}{2\pi } \int_{ \frac{\delta}{\sigma_{\max}} \leq |\xi| } \exp \left( - \frac{\sigma^2 |\xi|^2}{2} \right) d \xi  \\
  & \leq \frac{1}{2\pi}  \int_{\frac{\delta}{\sigma_{\max}} \leq |\xi|  }  \left| \mathbb{E}\left[ \exp \left(i  \sum_{k=1}^n  (kX_k - \mathbb{E} \left[k X_k\right])  \xi \right) \right]\right|  d \xi  +  \frac{1}{2\pi } \int_{ \frac{\delta}{\sigma_{\max}} \leq |\xi| } \exp \left( - \frac{\sigma^2 |\xi|^2}{2} \right) d \xi  \\
& =: T_3 + T_4 \label{e_estiamte_T_2_triangle}
\end{align}
We will estimate each term on the right hand side of the last estimate separately. We start with estimating the term~$T_4$. Direct calculation yields
\begin{align}
  T_4 & \leq C \frac{1}{\sigma} \int_{ |\xi| \geq \frac{\delta}{\sigma_{\max}} \sigma}  \exp \left( - \frac{ |\xi|^2}{2} \right) d \xi 
 \leq C \frac{1}{\sigma} \, \left(\frac{\sigma_{\max}}{\delta \sigma}\right) \exp \left( - \frac{\sigma^2\delta^2}{2\sigma_{\max}^2} \right)
 \leq C \frac{\sigma_{\max}}{\sigma^2}. \label{e_estimate_T_4}
\end{align}

Now, let's turn to the estimation of the term~$T_3$ on the right hand side of~\eqref{e_estiamte_T_2_triangle}. Using the trivial estimate that for~$1 \leq l \leq n $
\begin{align*}
\left|   \mathbb{E} \left[ \exp \left( i (lX_l - \mathbb{E} [lX_l]) \xi  \right)  \right] \right| \leq 1
\end{align*}
we get that
\begin{align}
T_3 & = \frac{1}{2\pi}  \int_{\frac{\delta}{\sigma_{\max}} \leq |\xi|  }  \left| \mathbb{E}\left[ \exp \left(i  \sum_{k=1}^n  (kX_k - \mathbb{E} \left[k X_k\right])  \xi \right) \right]\right|  d \xi  \\
  & \leq \frac{1}{2\pi}  \int_{\frac{\delta}{\sigma_{\max}} \leq |\xi|  }  \left| \mathbb{E}\left[ \exp \left(i  \sum_{k\in M}  (kX_k - \mathbb{E} \left[ kX_k\right])  \xi \right) \right]\right|  d \xi \\
\label{eq:T3} & = \frac{1}{2\pi}  \int_{\frac{\delta}{\sigma_{\max}} \leq |\xi|  }  \left| \mathbb{E}\left[ \exp \left(i  \sum_{k\in M}  \frac{C_3}{C_{\xi,k} \sqrt{\var{X_k}}}  {(X_k - \mathbb{E} \left[ X_k\right])} \  \frac{k \sqrt{\var{X_k}}}{\sigma_{\max}}   \frac{C_{\xi,k}}{C_3}\sigma_{\max} \xi \right) \right]\right|  d \xi, 
\end{align}
where the constant~$0< C_3 < \infty$ is given by~\eqref{e_def_set_M} and~$0<C_{\xi,k}< \infty$ is defined now. 
Recall that $\sigma_k^2 = \var{k\, X_k}$. 
We introduce the function~$|\cdot|_{2\pi}: \mathbb{R} \to [- \pi, \pi]$ via
\begin{align*}
  |x|_{2 \pi} = |r| , \qquad \mbox{where~$r$ is determined by } ~x= k 2 \pi + r \quad \mbox{with } k \in \mathbb{Z}.
\end{align*}
Let $\tilde \delta>0$ be some fixed constant to be determined later. 
 For the definition of the constant~$C_{\xi,k}$, for a given~$\xi \in \mathbb{R}$ and~$k \in \mathbb{N}$, we distinguish two cases. In the first case we assume that
\begin{align*}
  \left| \frac{\sigma_k}{\sigma_{\max}} \frac{1}{C_3}  \sigma_{\max} \xi \right|_{2 \pi} \geq \tilde\delta
\end{align*}
and set
\begin{align*}
  C_{\xi,k}=1.
\end{align*}
In the case that
\begin{align}\label{e_def_C_xi_second_case}
 \left| \frac{\sigma_k}{\sigma_{\max}} \frac{1}{C_3} \sigma_{\max} \xi \right|_{2 \pi} < \tilde\delta
\end{align}
we define the constant~$C_{\xi,k}$ in the following way. Because
\begin{align*}
  |\xi| \geq \frac{\delta}{\sigma_{\max}} 
\end{align*}
and by definition (see~\eqref{e_def_set_M})
\begin{align*}
  \frac{1}{C_2} \leq \left| \frac{\sigma_k}{\sigma_{\max}} \right| \leq C_2
\end{align*}
It holds that
\begin{align*}
\left| \frac{\sigma_k}{\sigma_{\max}} \frac{1}{C_3} \sigma_{\max} \xi \right| \geq \frac{1}{C_3C_2} \delta \geq \tilde \delta,
\end{align*}
by choosing~$\tilde \delta>0$ small enough.
Together with~\eqref{e_def_C_xi_second_case} it follows that there is~$0 \neq m \in \mathbb{N}$ such that 
\begin{align*}
\left| \frac{\sigma_k}{\sigma_{\max}} \frac{1}{C_3} \sigma_{\max} \xi \right| = m 2 \pi + r. 
\end{align*}
Now we set
\begin{align*}
  C_{\xi,k} = 1 + \frac{1}{ m 2 \pi }.
\end{align*}
As a consequence we get that
\begin{align}\label{e_unform_bounds_xi_sigma_max}
  \left| C_{\xi,k} \frac{\sigma_k}{\sigma_{\max}} \frac{1}{C_3} \sigma_{\max} \xi  \right|_{2 \pi}  \geq 0.5.
\end{align}
Additionally, it follows from the definition of~$C_{\xi,k}$ that
\begin{align}\label{e_uniform_bounds_c_xi}
  1 \leq   \left| C_{\xi,k}\right| \leq 2 .
\end{align}
With these definitions, we can rewrite Equation~\eqref{eq:T3} as 
\begin{align*}
& \mathbb{E}\left[ \exp \left(i  \sum_{k\in M}  \frac{C_3}{C_{\xi,k} \sqrt{\var{X_k}}}  {(X_k - \mathbb{E} \left[ X_k\right])} \  \frac{k \sqrt{\var{X_k}}}{\sigma_{\max}}   \frac{C_{\xi,k}}{C_3}\sigma_{\max} \xi \right) \right]  
= \Pi_{k \in M} \mathbb{E} \left[  \exp \left(i Z_k  \zeta_k \right) \right],
\end{align*}
where the random variable~$Z_k$ is given by
\begin{align*}
 \frac{C_3}{C_{\xi,k} \sqrt{\var{X_k}}}  {(X_k - \mathbb{E} \left[ X_k\right])}
\end{align*}
and the number~$\zeta_k \in \mathbb{R}$ is given by
\begin{align*}
\zeta_k =  \frac{k \sqrt{\var{X_k}}}{\sigma_{\max}}   \frac{C_{\xi,k}}{C_3}\sigma_{\max} \xi .
\end{align*}
With the properties we have derived, one sees that~$Z_k$ satisfies the conditions of Lemma~\ref{p_control_characteristic_fnction_intermediate_and_large_values} uniformly in~$k$. Additionally, one can also verify using the properties from above that
\begin{align*}
|\zeta_k|_{2 \pi} \geq \tilde \delta.
\end{align*}
Therfore we can apply Lemma~\ref{p_control_characteristic_fnction_intermediate_and_large_values} and get that there is a constant~$0 < \lambda< 1 $ such that for all~$k \in M$
\begin{align*}
\left|  \mathbb{E} \left[ \exp \left( i Z_k  \zeta_k  \right)  \right] \right| \leq \lambda.
\end{align*}
Hence, overall we have 
\begin{align}
  T_3 & \leq \frac{1}{2\pi} \int_{\frac{\delta}{d_n} \leq |\xi|  }  \prod_{k=1}^n \left| \mathbb{E}\left[ \exp \left(i ( X_k - \mathbb{E}[X_k])  \sigma_k  \xi \right) \right] \right| d \xi \\ 
  \\& \leq \frac{1}{2\pi} \int_{\frac{\delta}{d_n} \leq |\xi|  }  \prod_{k=1}^n \left| \mathbb{E}\left[ \exp \left(i ( Z_k \zeta_k  \xi \right) \right] \right| d \xi \\
 & \leq \frac{1}{2 \pi} \lambda^{|M|}.\label{e_estimation_T_3_second_part}
\end{align}
Finally, a combination of the estimate~\eqref{e_estimate_T_4} and~\eqref{e_estimation_T_3_second_part} yields that for~$c = - \ln \lambda$, we have 
\begin{align*}
  T_2 \leq T_3 + T_4 \leq C \frac{\sigma_{\max}}{\sigma^2} + C e^{-c\, |M|} \,
\end{align*}
which closes the argument.

Finally, observe that when $h > 1$ is fixed, one simply replaces $\pi$ with $\pi/h$, and the corresponding analysis is the same. 
\end{proof}

\subsection{Quantitative Arratia--Tavare principle}
\label{AT:proof}

\begin{proof}[Proof of Theorem~\ref{AT:normal}]
Our proof is similar to the specializations in~\cite{PittelShape, PittelSetPartitions} for integer partitions and set partitions, in that we exploit Equation~\eqref{dtv:intro} by utilizing local central limit theorems; however, in those cases, the generating function was explicitly analyzed to obtain the LCLT, whereas we appeal to Theorem~\ref{p_qlclt_for_combinatorics} directly. 

Let $B \subset [n]$ denote some subset of indices.  
Let $T_1^\ast := (T_B - \mu_B)/\sigma_B$.  
Using Markov's inequality for any $t>0$, and Lemma~\ref{p_taylor_characteristic_function}, we have 
\begin{equation}\label{eq:ld} \P(T_B > n) \leq \exp\left(- \frac{n-\mu_B}{\sigma_B} \right) \e e^{t T_1^\ast} \leq \exp\left(- \frac{n-\mu_B}{\sigma_B} \right) \e \exp\left(\frac{t^2}{2} + O\left(|t|^3 \e |T_1^\ast|^3\right)\right).  \end{equation}
Equation~\eqref{eq:ld} motivates the first assumption, which is that $(n-\mu_B)/\sigma_B \to 0$. 
At this point we also require that $T_B$ satisfy a \emph{central limit theorem}, or sufficiently, a LCLT, so that the right-hand side of Equation~\eqref{eq:ld} is controlled. 
 Thus, there is a $C>0$ such that $\e |T_1^\ast|^3 \leq C < \infty$; taking $t=1$, we obtain a bound on the big-$O$ contribution, and so we have 
\[ \P(T_B > n) \leq C \exp\left(- \frac{n-\mu_B}{\sigma_B} \right).  \]

\ignore{Next, we assume
\begin{align}
 \mbox{ $\RAc$ satisfies a \emph{local} central limit theorem; } \\
 \mbox{ $T$ satisfies a \emph{local} central limit theorem. }
\end{align}}
Then, by Theorem~\ref{p_qlclt_for_combinatorics}, we have for all sets $B \subset [n]$, 
\begin{align}
 \P(T_B^c = n-r) & = \P\left(\frac{T_B^c-\mu_B^c}{\sigma_B^c} = \frac{n-\mu_B^c - r}{\sigma_B^c}\right) \\
  & = \frac{1}{\sigma_{B}^c \sqrt{2\pi}}\, \exp\left(-\frac{1}{2} \left(\frac{n-\mu_B^c - r}{\sigma_B^c}\right)^2  \right)\left(1 + O\left( \frac{\sigma_{B^c,\max}}{\sigma_{B^c}} + \sigma_{B^c}\, e^{-C_3 |M_B^c|}\right)\right)
\end{align}
and
\begin{align}
 \P(T = n) & = \P\left(\frac{T-\mu}{\sigma} = \frac{n-\mu}{\sigma}\right) \\
  & = \frac{1}{\sigma \sqrt{2\pi}}\, \exp\left(-\frac{1}{2} \left(\frac{n-\mu}{\sigma}\right)^2  \right)\left(1 + O\left( \frac{\sigma_{\max}}{\sigma} + \sigma\, e^{-C_3 |M|}\right)\right).
\end{align}
Whence, 
\begin{align}
\frac{\P(T_B^c=n-r)}{\P(T=n)} & = \frac{\sigma}{\sigma_B^c}\, \exp\left(-\frac{1}{2}\left(\frac{r-\mu_B+\mu_{err}}{\sigma_B^c}\right)^2 + \frac{\mu_{err}^2}{2\sigma^2} \right)\biggl(1 + \\
 & \qquad \qquad O\left( \frac{\sigma_{\max}}{\sigma} + \sigma\, e^{-C_3 |M|} + \frac{\sigma_{B^c, \max}}{\sigma_{B^c}} + \sigma_{B^c}\, e^{-C_3 |M_{B^c}|}\right)\biggr) \\
& = \frac{\sigma}{\sigma_B^c}\, \exp\left(-\frac{1}{2}\left(\frac{r-\mu_B}{\sigma_B}\right)^2\frac{(\sigma_B)^2}{(\sigma_B^c)^2} + \frac{\mu_{err}(r-\mu_B)}{(\sigma_B^c)^2}+ \frac{\mu_{err}^2}{2(\sigma_B^c)^2}+ \frac{\mu_{err}^2}{2\sigma^2} \right) \times \\
& \qquad\left(1 +  O\left( \frac{\sigma_{\max}}{\sigma} + \sigma\, e^{-C_3 |M|} + \frac{\sigma_{B^c, \max}}{\sigma_{B^c}} + \sigma_{B^c}\, e^{-C_3 |M_{B^c}|}\right)\right). 
\end{align}
We then have 
\begin{align}
\left| \frac{\P(T_B^c = n-r)}{\P(T=n)} - 1 \right| & = \left| \frac{\sigma}{\sigma_{B}^c} \exp\left(-\frac{(r-\e T_B)^2}{2(\sigma_{B}^c)^2}\right) - 1 \right| \times \biggl(1 + O\biggl(\frac{\mu_{err}|r-\mu_B|}{(\sigma_B^c)^2} + \frac{\mu_{err}^2}{(\sigma_B^c)^2} \\
& \qquad   + \frac{\sigma_{\max}}{\sigma} + \sigma\, e^{-C_3 |M|} + \frac{\sigma_{B^c, \max}}{\sigma_{B^c}} + \sigma_{B^c}\, e^{-C_3 |M_{B^c}|}\biggr)\biggr). 
\end{align}
Recalling that $T_1^\ast = (T_B - \mu_B)/\sigma_B$, we obtain 
\begin{align}
\label{dtv:preasymptotic} d_{TV}(C_B, Z_B) & = \frac{1}{2} \e \left| \frac{\sigma}{\sigma_B^c} \exp\left(-\frac{(T_1^\ast)^2}{2} \left(\frac{\sigma^2}{(\sigma_B^c)^2}\right) - 1 \right) - 1 \right| \times  \biggl(1+  \\
 & \qquad \qquad O\biggl( \exp\left(- \frac{n-\mu_B}{\sigma_B} \right) + \frac{\sigma_B\ \mu_{err}}{(\sigma_B^c)^2} + \frac{\mu_{err}^2}{(\sigma_B^c)^2}  \\
 & \qquad  \qquad +  \frac{\sigma_{\max}}{\sigma} + \sigma\, e^{-C_3 |M|} + \frac{\sigma_{B^c, \max}}{\sigma_{B^c}} + \sigma_{B^c}\, e^{-C_3 |M_{B^c}|}\biggr)\biggr). 
 \end{align}
 
Then under the assumptions
\begin{itemize}
\item $\displaystyle \sigma_{B^c}\, e^{-C_3 |M_B^c|}, \ \ \sigma\, e^{-C_3 |M|}, \ \ \frac{\sigma_{\max}}{\sigma}, \ \ \frac{\sigma_{B^c,\max}}{\sigma_{B^c}} \to 0$; 
\item $\displaystyle \mu_{err} = o(\sigma_B^c)$; \\
\item $\displaystyle \frac{\sigma_B \mu_{err}}{(\sigma_B^c)^2} \to 0$;\\
\item $\lim_{n\to\infty} \sigma_1 \in (0,1)$; 
\end{itemize}
and since $T_1^\ast$ is assumed to tend to a normal distribution, say at rate $r_B = o(1)$, we have 
\begin{align}
 d_{TV}(C_B, Z_B) & \to \frac{1}{2} \int_{-\infty}^\infty \frac{1}{\sqrt{2\pi}} \left| \frac{1}{\sigma_1} \exp\left(-\frac{x^2}{2} \left(\frac{1}{(\sigma_1)^2}\right) - 1 \right) - 1 \right| \exp\left(-\frac{x^2}{2}\right) dx \\
 & = \frac{1}{2} \int_{-\infty}^\infty  \left| \frac{1}{\sqrt{2\pi \sigma_1^2}} \exp\left(-\frac{x^2}{2\sigma_1^2}\right) - \frac{1}{\sqrt{2\pi}} \exp\left(-\frac{x^2}{2}\right) \right| dx \\
 & = d_{TV}(\mathcal{N}_B, \mathcal{N}),
 \end{align}
where $\mathcal{N}_B$ is a normal distribution with mean 0 and variance $\sigma_1^2$. 

When $\sigma_1\to 1$, then by Equation~\eqref{dtv:preasymptotic}, assuming $T_1^\ast$ has finite second moment, the expression in the right-hand side of Equation~\eqref{dtv:preasymptotic} tends to zero at a rate proportional to $\sigma_1^2 - 1$. 
\end{proof}

\section{Integer partitions}
\label{section:integer:partitions}
\subsection{Asymptotic enumeration}
\label{section:asymptotic:enumeration}
Let us first demonstrate various consequences of Theorem~\ref{p_qlclt_for_combinatorics} to some classical enumerative integer partition formulas. 
We start with unrestricted integer partitions, and then describe how to generalize the probabilistic model to reclaim many classical enumerative results for integer partitions under various restrictions, traditionally obtained by detailed analysis of the relevant generating functions. 

Taking $g_i(k) = 1$ in Equation~\eqref{eq:Z}, for any $0<x<1$ we have 
\[ \P(Z_i = k) = c_i\, x^{i\,k} = x^{i\, k} (1-x^i), \qquad i=1,2,\ldots, \quad  k=0,1,\ldots,\]
hence $Z_i$ is a geometric random variable with parameter $1-x^i$, $i=1,2,\ldots$.  
Interestingly, \emph{for all $0<x<1$}, we also have (see for example~\cite{Fristedt, Vershik})
\begin{equation}\label{ip:pn} \P\left(\sum_{i=1}^n i\, Z_i = n\right) = p(n) x^n \prod_{i=1}^n (1-x^i),  \qquad n \geq 0, \end{equation}
where $p(n)$ denotes the number of integer partitions of size~$n$. 
From this one obtains the interpretation that conditional on the event $\{\sum_{i=1}^n i\, Z_i = n\}$, $Z_i$ is the number of parts of size~$i$ in a uniformly random integer partition of $n$. 

Let $c := \pi/\sqrt{6}$.  
We shall take $x = \exp(-c/\sqrt{n})$ in what follows, which is a highly strategic choice. 
We have 
\[ \e T = \sum_{i=1}^n \frac{i\, x^i}{1-x^i} = \sum_{i=1}^n \frac{i\, e^{-c\,i/\sqrt{n}}}{1-e^{-c\,i/\sqrt{n}}}. \]
Using Riemann sums, with $y_i = i/\sqrt{n}$ and $\Delta y_i = 1/\sqrt{n}$, we have 
\[ \e T - n = O(\sqrt{n}). \]
We also have 
\[ \var{i\, Z_i} = \frac{i^2 x^i}{(1-x^i)^2} = \frac{i^2 e^{-c\,i/\sqrt{n}}}{(1-e^{-c\,i/\sqrt{n}})^2}, \]
and so 
\[ \var{T} =  \sum_{i=1}^n \frac{i^2\, e^{-c\,i/\sqrt{n}}}{(1-e^{-c\,i/\sqrt{n}})^2}. \]
Again appealing to Riemann sums, we have 
\[ \sigma^2 := \var{T} \sim \frac{2}{c}\, n^{3/2}. \]
In order to apply Theorem~\ref{p_qlclt_for_combinatorics}, we also need $\sigma_{\max}$ and to find a set $M$.  
First, for $i = o(\sqrt{n})$, we have 
\[ \var{i\, Z_i} = \frac{i^2 e^{-c\,i/\sqrt{n}}}{(1-e^{-c\,i/\sqrt{n}})^2} \sim \frac{n}{c^2}. \]
For $i \sim t\, \sqrt{n}$, $t > 0$, we have 
\[ \var{i\, Z_i} = \frac{i^2 e^{-c\,i/\sqrt{n}}}{(1-e^{-c\,i/\sqrt{n}})^2} \sim \frac{n}{c^2}\frac{t^2 e^{-c\,t}}{(1-e^{-c\,t})^2}, \]
where $\frac{t^2 e^{-c\,t}}{(1-e^{-c\,t})^2}$ is decreasing for increasing $t>0$. 
Finally, for $i/\sqrt{n} \to \infty$, we have $\var{i\, Z_i} = o(n)$, and so $\sigma_{\max}/\sigma = O(\sqrt{n}/n^{3/4}) = O(1/n^{1/4})$. 
In addition, for some $C>0$, we have $M = \{1, 2, \ldots, k_1\}$, where $k_1 = \Omega(\sqrt{n})$; i.e., $|M| = \Omega(\sqrt{n})$.
Applying Theorem~\ref{p_qlclt_for_combinatorics}, we have 
\begin{equation}\label{almost:HR} p(n) e^{-c\, \sqrt{n}} \prod_{i=1}^n \left(1-e^{-c\, i/\sqrt{n}}\right) =  \P(T = n) = \frac{1}{\sqrt{2\pi\, \frac{2}{c}\, n^{3/2}}}\left(1 + O\left(\frac{1}{n^{1/4}} + n^{3/4} e^{-C_1 \sqrt{n}}\right)\right). \end{equation}
We are almost finished. 
The last step is to expand out the product on the left-hand side of Equation~\eqref{almost:HR}. 
This type of analysis can be found in, e.g.,~\cite[Chapter~2]{newman2006analytic}, \cite{Ingham, Romik, NonNegR}, via a variety of methods. 
In the unrestricted case above, this product is~\cite[Exercise~3, page~58]{Debruijn}. 
 Rearranging, this gives us an approximation for $p(n)$, namely, 
\begin{equation}\label{eq:HR} p(n) = \frac{e^{2c\sqrt{n}}}{4\sqrt{3}\, n}\left(1 + O\left(\frac{1}{{n^{1/4}}}\right)\right). \end{equation}

The leading term of the right-hand side above is sometimes referred to as the Hardy--Ramanujan asymptotic formula for $p(n)$, see~\cite{HR}, although their results are considerably stronger; see for example~\cite[Section~7]{logconcave} for a more thorough discussion. 

Also of classical interest in the enumeration of combinatorial sequences is the number of integer partitions of size~$n$ into exactly $k$ parts. 
Fortunately, the number of these partitions is in bijection with the number of integer partitions of size~$n$ in which the largest part size is~$k$, which we denote by $p(n,k)$. 
More generally, we could let $U = \{u_1, u_2, \ldots \} \subset \{1,2,\ldots\}$ denote a subset of allowable part sizes, and denote by $p_U(n)$ the number of integer partitions of $n$ into parts with sizes in the set $U$. 
This generalization has been studied extensively, see for example~\cite{nathanson2000partitions, Szekeres1, Szekeres2, RothSzekeres, erdos1976concerning, ErdosElementary, Wright, Ingham, hua1942number, Freiman, Romik, Meinardus}.  
We can generalize even further by restricting the multiplicities of the part sizes, see for example~\cite{hagis1971partitions, CanfieldWilf}; that is, we similarly let $R \subset \{0,1,2,\ldots\}$ denote a set of allowable part sizes (taking $R = \{0,1\}$, for example, corresponds to partitions into distinct parts). 
Each of these generalizations and specializations is obtained by a careful analysis of the relevant generating functions, with formulas containing similar types of expressions.  
The probabilistic model offers a universal approach in many settings. 

There are several regimes of interest, and we first recount some of the main results from the literature.  
First, let us assume $R = \{0, 1, 2, \ldots \}$, i.e., no restrictions on part sizes. 

\begin{theorem}[see, e.g.,~\cite{nathanson2000partitions}]\label{finite:theorem}
Suppose $U \subset \{1, 2, \ldots \}$ is a finite set of precisely $k$ elements which are relatively prime. 
Then as $n$ tends to infinity we have  
\[ p_U(n) = \frac{1}{\prod_{a\in U} a} \frac{n^{k-1}}{(k-1)!} + O\left(n^{k-2}\right). \]
\end{theorem}
A similar asymptotic formula also holds for $p(n,k)$. 
\begin{theorem}[\cite{ErdosLehner}]
Let $U = \{1, 2, \ldots, k\}$ for some positive integer $k$.  Then as $n$ tends to infinity we have 
\[ p(n,k) \sim \frac{1}{k!}\binom{n-1}{k-1} \]
uniformly in $k$ for $k = o(n^{1/3})$.
\end{theorem}
This was further extended by Szekeres.
\begin{theorem}[\cite{Szekeres1}]
Let $U = \{1, 2, \ldots, k\}$ for some positive integer $k$.  Then as $n$ tends to infinity we have 
\[ p(n,k) \sim \frac{1}{k!}\binom{n-1}{k-1} \exp\left(\frac{1}{4}\, \frac{k^3}{n} - \frac{13}{288} \frac{k^5}{n^2} + \cdots\right) \]
uniformly in $k$ for $k = o(n^{1/2})$.
\end{theorem}
These formulas indicate sub-exponential growth if we restrict the allowable part sizes too much. 
It was shown in~\cite{ErdosLehner} that the largest part in a random unrestricted integer partition is on the order of $\frac{1}{2c} \sqrt{n} \log(n)$. 
Thus, taking $k$ larger than this value, we have $p(n,k) \sim p(n)$. 
However, if we only include all part sizes up to $t \sqrt{n}$, for some $t>0$, then we have the following. 
\begin{theorem}[\cite{Romik}]
For any $t>0$, let $U_t = \{1, 2, \ldots, [t\sqrt{n}]\}$.  then we have 
\begin{equation}\label{eq:Romik} p_{U_t}(n) \sim \frac{G(t)}{n} \, e^{H(t)\sqrt{n}}, \end{equation}
where 
\[ G(t) = \frac{\alpha(t)}{2\pi\left(2-(t^2+2)e^{-\alpha(t)t}\right)^{1/2}}, \]
\[ H(t) = 2\alpha(t) - t \log(1-e^{-\alpha(t) t}), \]
and $\alpha(t)$ satisfies 
\[ \alpha(t)^2 = \Li_2(1-e^{-\alpha(t)t}), \]
where 
\[ \Li_2(x) := \sum_{m=1}^\infty \frac{x^m}{m^2} \]
is the dilogarithm function. 
\end{theorem}

One can also ask for partitions of $n$ with all part sizes exceeding a given value $r$; see for example~\cite{almkvist2002partitions, dixmier1990partitions, dixmier1990partitions2, nicolas2000partitions}.  
When $r=1$, there is a direct argument as follows: the number of partitions of $n$ with no parts of size~$1$ is equal to $p(n) - p(n-1)$, and so for $U = \{2, 3, \ldots \}$, we have
\begin{equation}\label{eq:no:one} p_U(n) = p(n)-p(n-1) \sim p(n) \frac{c}{\sqrt{n}}.
\end{equation}
However, simple arguments of this form do not always generalize so simply. 
We have chosen to present the following special case of a more general theorem in~\cite{nicolas2000partitions}.
We refer the interested reader to~\cite{nicolas2000partitions} for the more comprehensive treatment, which provides higher order terms and a more detailed treatment of the function $\tilde g$ in the theorem below. 
\begin{theorem}[\cite{nicolas2000partitions}]\label{theorem:nicolas}
Suppose $U = \{r, r+1, \ldots \}$.  Then we have for $m \leq \sqrt{n}$, 
\begin{equation}\label{eq:Nicolas} p_U(n) \sim p(n) \left(\frac{2c}{\sqrt{n}}\right)^{r-1} (r-1)!\, \exp\left(\sqrt{n}\, \tilde g\left(\frac{r^2}{\sqrt{n}}\right) \right), \end{equation}
where $\tilde g$ is made explicit in several regions for $m$. 
\end{theorem}

We can also restrict to subsets of the integers with positive density strictly less than 1, for example partitions into parts of sizes in set $J$ modulo $m$. 
Taking $J=\{1\}$ and $m=2$, we obtain partitions into odd parts; taking $J = \{0\}$ and $m=2$, we obtain partitions into even parts. 
Many partition bijections are given in terms of sets of partitions that fit into this framework, see~\cite{Andrews, PakSurvey}. 

\begin{theorem}[\cite{robertson1976partitions}]
Suppose positive integer $m>0$ is given, along with set $J = \{j_1, \ldots, j_k\}$, with $k \leq m$.  
Let $U = \{J + \ell\, m: \ell=0,1,2,\ldots\}$ be such that the elements $u_k$ modulo $m$ lie in the set $J$. 
Then we have 
\[ p_U(n,k) \sim \frac{n^{k-1}}{k!(k-1)!} \frac{|J|^k \delta}{m}, \qquad n \to \infty, \quad k = o(n^{1/4}),\]
where $\delta$ is the greatest common divisor of $(j_2-j_1, \ldots, j_k - j_1, m)$. 
\end{theorem}

\begin{theorem}[\cite{Meinardus}]
Suppose relatively prime positive integers $k$ and $a$ are given, with $1 \leq a < k$. 
Let $U = \{a + \ell\, k: \ell=0,1,2,\ldots\}$ be such that the elements $u_k \equiv a$ modulo $m$. 
Then we have 
\[ p_U(n) \sim C n^{\kappa} e^{\pi \sqrt{2n/3k}}, \]
where
\[ C = \Gamma\left(\frac{a}{k}\right)\cdot  \pi^{\frac{a}{k}-1}\cdot 2^{-\frac{3}{2}-\frac{a}{2k}}\cdot 3^{-\frac{a}{2k}}\cdot k^{-\frac{1}{2}+\frac{a}{2k}} \]
and
\[ \kappa = -\frac{1}{2}\left(1+\frac{a}{k}\right). \]
\end{theorem}

A second popular generalization considers sets $U$ for which $u_k$ grows polynomially; roughly, $u_k \sim B k^r$ for some $B>0$ and $r \geq 1$; we say the sequence $U$ has \emph{polynomial growth} in this case, although care must be taken since the higher order asymptotic terms of $u_k$ affect the final asymptotic formula for $p_U(n)$. 
There is also quite a history in this setting as well, starting with Hardy and Ramanujan~\cite{HR} who suggested their method could be generalized for $B=1$ and $r$ a positive integer. 
The work in~\cite{Wright} carries out this analysis, and in~\cite{Ingham} a Tauberian theorem is proved and applied to integer partitions satisfying a more refined polynomial growth condition. 

\begin{theorem}[\cite{Ingham}]\label{ingham:theorem}
Suppose $U = \{u_1, u_2, \ldots \}$ satisfies the following. 
Let $N(u)$ denote the number of elements in $U$ less than or equal to $u$.  Suppose there are constants $L>0$, $\theta >0$, and function $R(u)$ such that 
\[ N(u) = L u^\theta + R(u), \]
and for some $b>0$ and $c>0$ we also have 
\[ \int_0^u \frac{R(v)}{v}\, dv = b \log u + c + o(1). \]
Define
\[ \alpha = \frac{\beta}{1+\beta}\qquad M = \{L \beta\Gamma(\beta+1)\zeta(\beta+1)\}^{1/\beta}. \]
Then as long as $\gcd(U) = 1$ (a simpler sufficient condition provided later by~\cite{ErdosBateman}, also observed in~\cite{NonNegR}), we have 
\begin{equation}\label{eq:Ingham} p_U(n) \sim \left(\frac{1-\alpha}{2\pi}\right)^{1/2}e^cM^{-(b-\frac{1}{2})\alpha}u^{(b-\frac{1}{2})(1-\alpha)-\frac{1}{2}}e^{\alpha^{-1}(Mu)^\alpha}. \end{equation}
\end{theorem}

Further developments are contained in~\cite{RothSzekeres}, where the conditions in Theorem~\ref{ingham:theorem} are replaced by 
\begin{equation}\tag{U1} \lim_{k\to \infty} \frac{\log u_k}{\log k} = r \ \ \mbox{ exists,} \end{equation}
\begin{equation}\tag{U2} 
J_k^\ast = \inf_{\frac{1}{2}u_k^{-1} < \alpha \leq \frac{1}{2}}\left\{(\log k)^{-1}\sum_{\nu=1}^k \left(\frac{u_k}{u_\nu}\right)^2 [u_\nu \alpha]^2 \right\} \to \infty.
\end{equation}
These conditions allow for $u_k$ to be the $k$-th prime number, for example, and the formulas stated in~\cite[Section~5]{RothSzekeres} contain higher order terms of an asymptotic expansion.  

We can also allow certain cases where $u_k$ is allowed to grow exponentially fast. 
The problem appears to have originally been proposed and developed by Mahler~\cite{mahler1940special} for $u_k = \nu^k$, for $\nu >0$, with a later refinement by deBruijn~\cite{Debruijn}. 
In this case, Mahler showed
\[ \log p_U(n) \sim \frac{1}{2\log \nu} \left(\log \frac{n}{\log{n}}\right)^2 + \left(\frac{1}{2}+\frac{1}{\log\nu} + \frac{\log \log \nu}{\log \nu}\right) \log{n} - \left(1+\frac{\log \log \nu}{\log \nu}\right)\log\log{n} + O(1), \]
and deBruijn later refined the $O(1)$ estimate as being of the form 
\[ \psi\left( \frac{\log n - \log\log{n}}{\log \nu} \right) + o(1), \]
where $\psi$ is a certain periodic function of period 1.
Richmond~\cite{Richmond} also has results in this area, with a generalization and assumptions akin to~\cite{RothSzekeres}. 

\begin{theorem}[\cite{Richmond}] \label{theorem:Richmond}
Suppose $U = \{u_1, u_2, \ldots \}$ is given. 
Define $A(u)$ for $u>0$ as the number of elements of $U$ which are less than or equal to $u$. 
Assume 
\begin{itemize}
\item[(i)] $ s = \lim_{k\to\infty} \frac{\log\log{u_k}}{\log{k}}$ exists;
\item[(ii)] $\liminf_{k\to\infty} \frac{\log{u_k}}{k} > 0$; 
\item[(iii)] $A(2k) = O(A(k))$ as $k \to \infty$;
\item[(iv)] $\gcd(U) = 1$;
\item[(v)] $\limsup \frac{\log\log u_k}{\log{k}} < \infty$.
\end{itemize}
When $U$ satisfies condition~$(i)$ we say that $U$ has \emph{exponential growth}.  
Using the notation from~\cite{Richmond}, let $\alpha$ denote the solution to 
\[ n = \sum_k u_k (e^{\alpha u_k}-1)^{-1}, \]
and let
\[ A_2 = \sum_k \frac{u_k e^{-\alpha u_k}}{(1-e^{-\alpha u_k})^2}. \]
Then for any integer $m \geq 2$, we have
\begin{align}\label{eq:Richmond}  
p_U(n) & =  \frac{1}{\sqrt{2\pi A_2}} \exp\left(\sum_{\nu=0}^\infty \left[\frac{\alpha u_k}{e^{\alpha u_k}-1} - \log\left(1-e^{-\alpha u_k}\right)\right]\right) \times  \\
 & \qquad \qquad  \left(1 + \sum_{\rho=1}^{m-2} D_\rho + O\left(f_U^{1-(2m/3)}(\alpha)\right) \right), 
\end{align}
where expressions for $D_\rho$ can be found in~\cite[p.~390]{Richmond}, and the function $f_U(x) = \sum_k e^{-x u_k}$ for $x>0$.

\end{theorem}

The expressions for $\alpha$ and $A_2$ should of course look familiar, as our strategic choice of $x$ earlier in the section is an approximation to $\alpha$, and $A_2$ is precisely the variance of the random variable $T$. 
In fact, this is a common theme that permeates throughout much of the asymptotic analysis of integer partitions. 

We can also restrict the multiplicity of part sizes; typically, this is in the form of partitions into distinct parts. 
Similar theorems hold in this case as well, as the two problems are usually treated together as complementary combinatorial structures.  
We shall not recount the similar theorems in this setting, since most of the previously referenced work also contains theorems applicable to partitions into distinct parts, but there are several noteworthy generalizations.  
In this respect, we define, similarly as in~\cite{CanfieldWilf} the quantity $p(n, U, R)$, which is the number of partitions of $n$ into part sizes from the set $U$ and multiplicities in the set $R$. 
\begin{theorem}[\cite{hagis1971partitions}]\label{theorem:Hagis}
Suppose $U = \{1, 2, \ldots, \}$, and $R = \{0, 1, 2, \ldots, t\}$, i.e., no part can be repeated more than $t$ times. 
Let $r = t/24$, $E = 4\pi(n+r)^{1/2}$, $s=(t+1)^{-1/2}$. 
Then we have 
\begin{equation}\label{eq:Hagis} p(n,U,R) = \frac{\sqrt{12}\, s^{3/2}\, t^{1/4}}{(24n+t)^{3/4}}\, \exp\left(\frac{s\, E}{\sqrt{24}} \sqrt{n}\right)(1+O(n^{-1/2})). \end{equation}
\end{theorem}

In particular, taking $t=1$ in Theorem~\ref{theorem:Hagis}, we obtain integer partitions into distinct part-sizes. 
Our probabilistic model is well-equipped to handle restrictions of the form involving sets $U$ and $R$; that is, partitions where part sizes must come from set $U$ and can only have multiplicities in set $R$; in fact, we can further generalize so that each element $i \in U$ has its own set of allowable multiplicities $R_i$.  
The reason for the ease in generalization lies in the properties of the geometric random variables composing the sum $T$. 
For example, let $Z$ denote a geometric random variable. 
For any positive integer $r$, we have 
\[ \P(Z = k | Z \leq r) = \frac{x^k}{1+x+x^2 + \ldots, x^r}, \]
where in the special case of $r=1$ we have $\L(Z | Z \leq 1) =$ Binomial$(x/(1+x))$. 
In the other direction, we have $\L(Z | Z \geq r) = \L(r+Z)$; that is, the distribution of $Z$ conditioned on being at least $r$ is equal to the distribution of $Z$ shifted up by $r$. 

In general, given sets $U$ and $R_1, R_2, \ldots$ of integers, Equation~\eqref{eq:Z} then becomes 
\[ \P(Z_i = k) = c_i\, x^{i\,k}\, \mathbbm{1}(i \in U, k \in R_i), \qquad i\in U, \ k=0,1,2,\ldots, \]
and $\P(Z_i = 0) = 1$ for $i \notin U$. 
In other words, $g_i(k) = \mathbbm{1}(i \in U, k \in R_i)$, where of course $R_i$ must contain $0$ whenever $i \notin U$. 
Generalizing Equation~\eqref{ip:pn} to arbitrary $R_i$, we have 
\begin{equation}\label{prob:general} \P\left(\sum_{i \in U} i\, Z_i = n\right) = p_U(n) x^n \prod_{i\in U}\frac{1}{\sum_{j \in R_i} x^j}. \end{equation}

Thus, from a local central limit theorem for the sum $\sum_{i\in U} i\, Z_i$ we may estimate the quantity $p_U(n)$ by an analysis of the product in the right-hand side of Equation~\eqref{prob:general}. 
We shall not attempt such an analysis in full generality at present. 

First, however, is the admission that our methods require that the cardinality of $M$ must be tending to infinity in $n$, and so we must exclude at present those enumerative results pertaining to partitions into a finite subset of part sizes, i.e., Theorem~\ref{finite:theorem}. 
Furthermore, we must also check that $\sigma_{\max} / \sigma \to 0$ as well as $\sigma e^{- C |M|} \to 0$. 
For unrestricted integer partitions, it is sufficient to take $M = \{1, 2, \ldots, n^a\}$ for any $0 < a \leq 1/2$, i.e., the set $M$ contains polynomially many elements. 

\begin{theorem}\label{integer:partition:theorem}
Suppose $U = \{u_1, u_2, \ldots \}$ is a given subset of $\{1, 2, \ldots\}$ with $\gcd(U) = 1$. 
For each $i \in U$ suppose $R_i = \{0,1,\ldots\}$ or $R_{i} = \{0, 1, \ldots r\}$ for some positive integer $r$. 
Let ${\bf R} := (R_1, R_2, \ldots)$. 
For each $i=1,2,\ldots$, suppose $Z_i$ has distribution~\eqref{eq:Z} with $g_i(k) = \mathbbm{1}(i \in U, k \in R_i)$. 
Define $T(U,{\bf R}) := \sum_{i\in U} i\, Z_i$.
Suppose $x$ is such that 
\[ \e T(U,{\bf R}) - n =  \sum_{i \in U} \sum_{j \in R_i} j x^{i\,j} - n = o(\sigma), \]
where $\sigma^2 = \var{T(U,{\bf R})}$. 
Finally, suppose $M \subset U$ is some subset of indices which satisfies~\eqref{e_def_set_M} for a given constant $C>0$. 
Let $p(n; U,{\bf R})$ denote the number of integer partitions of $n$ into parts in the set $U$ with multiplicities restricted to lie in the set ${\bf R}$. 
Then asymptotically as $n$ tends to infinity, we have 
\[ p(n; U,{\bf R}) = \frac{x^{-n} \prod_{i\in U} \left(\sum_{j \in R_i} x^j\right)^{-1}}{\sqrt{2\pi \sigma^2}}\left(1 + O\left(\frac{\sigma_{\max}}{\sigma} + \sigma e^{-C |M|}\right) \right). \]
\end{theorem}
\begin{proof}
The main assumption is that the function $g_i(k)$ is perturbed log-concave, which is satisfied for indicator random variables over a contiguous set of indices. 
\end{proof}

\begin{remark}
The condition on $R_i$ can be replaced by any contiguous subset of indices by the following argument. 
For each $i\in U$, suppose $r_1^\ast$ is the first non-zero element of $R_i$: 
replace $R_i$ with $R_i - r_1^\ast$, and replace $n$ by $n-r_1^\ast$. 
Then, if still applicable, we may apply Theorem~\ref{integer:partition:theorem}. 
\end{remark}

\begin{remark}
A perhaps ``simpler" parameterization would be to forego the set $U$ altogether and allow the sets $R_i$ to contain just the element $\{0\}$. 
This notation would be very much at odds with classical treatments of integer partitions, which have focused primarily on restricting the set of part sizes and not their multiplicities. 
\end{remark}

\begin{remark}
We also strongly note that we are allowed to let the sets $U$ and ${\bf R}$ also grow with $n$, e.g., as in partitions of $n$ into parts of size at most $t \sqrt{n}$, i.e., $R_i = \{1,2,\ldots, t\sqrt{n}\}$. 
\end{remark}

We obtain a large number of corollaries, each of which corresponds to the first term in the asymptotic expansion of many previous results subject to an asymptotic analysis of the product $\prod_{i\in U}\sum_{j \in R_i}x^j$. 
In what follows, when all $R_i$ are taken to be the same, we refer to their common value as $R$.
Finally, many of the previously obtained results were obtained via arguments which also yielded higher order terms; in these cases, we only claim to agree with the first order asymptotic expression, with a valid big-$O$ error term. 
\begin{itemize}
\item $U = \mathbb{N}$, $R = \mathbb{N}\cup \{0\}$: see Equation~\eqref{eq:HR};
\item $U = \mathbb{N}$, for any $t>0$, $R = \{1, 2, \ldots, t\sqrt{n}\}$: see Equation~\eqref{eq:Romik};
\item $\log(u_k)/\log k \to r$ exists, $\gcd(U) = 1$, $R = \mathbb{N}\cup \{0\}$: see Equation~\eqref{eq:Ingham};
\item $U = \mathbb{N}$, for any $t>0$, $R = \{1, 2, \ldots, t\}$: see Equation~\eqref{eq:Hagis}. 
\end{itemize}

\begin{example}{\rm \label{example:no:small}
A slightly more involved corollary yields Equation~\eqref{eq:Nicolas}, which corresponds to the event  $\{Z_1 = Z_2 = \ldots Z_{r-1} = 0\}$. 
Let $p_r(n)$ denote the number of partitions of $n$ into parts of sizes at least $r$, and let $T_r := \sum_{i=r}^n i\, Z_i$.  
For any $0<x<1$, we have 
\begin{align} 
\P(Z_1 = Z_2 = \ldots Z_{r-1} = 0) & = \prod_{i=1}^{r-1} (1-x^i), \\
\P(Z_1 = Z_2 = \ldots Z_{r-1} = 0 | T=n) & = \frac{p_r(n)}{p(n)}, \\
\P(T=n, Z_1 = \ldots Z_{r-1}=0) & = p_r(n) x^n \prod_{i=1}^n (1-x^i).
\end{align}
Since
\[ \P(T = n | Z_1 = \ldots = Z_{r-1} = 0) = \P\left(\sum_{i=r}^n i\, Z_i = n\right) = p_r(n) x^n \prod_{i=r}^n (1-x^i), \]
and by Bayes' theorem we also have 
\[ \P(T = n | Z_1 = \ldots = Z_{r-1} = 0) = \frac{\P(Z_1 = \ldots Z_{r-1}=0 | T=n) \P(T=n)}{\P(Z_1 = Z_2 = \ldots Z_{r-1} = 0)} = \frac{p_r(n)}{p(n)} \ \frac{\P(T=n)}{\prod_{i=1}^{r-1} \P(Z_i = 0)}, \]
we apply Theorem~\ref{p_qlclt_for_combinatorics} for any $r = O(\sqrt{n})$, and obtain (with $x = e^{-c/\sqrt{n}}$)
\begin{align}
p_r(n) & = p(n) \frac{\P(T_r=n)}{\P(T=n)} \prod_{i=1}^{r-1}(1-x^i) = p(n)\ \left(\frac{c}{\sqrt{n}}\right)^{r-1} (r-1)!\ \sqrt{\frac{\var{T}}{\var{T_r}}}(1+O(n^{-1/4})) \\
&= p(n)\ \left(\frac{c}{\sqrt{n}}\right)^{r-1} (r-1)!\  \sqrt{\frac{2/c}{\int_r^\infty \frac{y^2 e^{-c\, y}}{(1-e^{-c\, y})^2}\,dy}}\left(1+O(n^{-1/4})\right)\ . \\
\end{align}
}\end{example}

\begin{example}{\rm \label{example:gap}
We shall also be interested in various statistics related to integer partitions in the next section, in particular the notion of a \emph{gap} in a partition. 
A partition is said to have a \emph{least gap} at $r$ if there are no parts of size~$r$ and at least one part each of sizes $1,2,\ldots, r-1$. 
A random partition with least gap at $r$ can be obtained by applying the transformations $Z_i \mapsto 1+Z_i$, $i=1,2,\ldots, r-1$, and $n \mapsto n-\binom{r}{2}$, and considering the event $\{Z_r = 0\}$. 
Letting $p_{\{r\}}(n)$ denote the number of partitions of $n$ with least gap equal to $r$, we have 
\begin{align}
\P(Z_1>0, \ldots Z_{r-1} >0, Z_r = 0) & = x^{\binom{r}{2}}\prod_{i=1}^{r} (1-x^i), \\
\P(Z_1 >0,  \ldots Z_{r-1} > 0, Z_r = 0 | T=n) & = \frac{p_{\{r\}}(n)}{p(n)}, 
\end{align}
\begin{align}\label{partitions:gap:pre}
\P\left(T=n, Z_1>0, \ldots Z_{r-1} >0, Z_r = 0\right) = p_{\{r\}}\left(n\right) x^{n} \prod_{i=1}^n(1-x^i). 
\end{align}
\begin{align}\label{partitions:gap}
\P\left(T=n| Z_1>0, \ldots Z_{r-1} >0, Z_r = 0\right) = \P\left(\sum_{i=1, i\neq r}^{n}i\, Z_i = n-\binom{r}{2}\right).
\end{align}
Applying Theorem~\ref{p_qlclt_for_combinatorics} to Equation~\eqref{partitions:gap}, for any $r < \sqrt{2n}$, we have 
\begin{align}
p_{\{r\}}(n) & = p(n) e^{-c\, \binom{r}{2}/\sqrt{n}} \left(\frac{c}{\sqrt{n}}\right)^r r!\, \sqrt{\frac{\var{T_n}}{\var{T_{n-\binom{r}{2}}}}}\left(1+O(n^{-1/4}) \right) \\
 & = p(n) e^{-c\, \binom{r}{2}/\sqrt{n}} \left(\frac{c}{\sqrt{n}}\right)^r r!\, \left(1+O(n^{-1/4}) \right). 
 \end{align}
Applying Theorem~\ref{p_qlclt_for_combinatorics} to Equation~\eqref{partitions:gap}, for $r$ such that $n - \binom{r}{2} \sim t \sqrt{n}$, for any $t>0$, we have 
\begin{align}
p_{\{r\}}(n) & = p(n) e^{-c\, \binom{r}{2}/\sqrt{n}} \left(\frac{c}{\sqrt{n}}\right)^r r!\, \sqrt{\frac{2/c}{\int_{0}^t \frac{y^2 e^{-c\, y}}{(1-e^{-c\, y})^2}\,dy}}\left(1+O(n^{-1/4})\right)\ . 
\end{align}
}\end{example}

\ignore{\begin{theorem}
Assume the conditions of Theorem~\ref{integer:partition:theorem}.  
Suppose for all $i$, $R_i = \{0, 1, 2, \ldots \}$, and $F_i$ is any finite subset of $\{1,2,\ldots\}$. 
Then Theorem~\ref{integer:partition:theorem} continues to hold with $R_i$ replaced with $R_i \setminus F_i$. 
\end{theorem}
\begin{proof}
Since all quantities involved are obtained via Riemann sums, modifying a finite number of elements has no impact on the asymptotic evaluations. 
\end{proof}}

\begin{remark}
The form of Theorem~\ref{theorem:Richmond} suggests that our analysis also applies in the case when $u_k$ grows exponentially. 
We have not carried out the detailed analysis, and indeed such analyses also appear to be absent from the limit shape literature. 
\end{remark}

The types of restrictions thus far encountered are referred to as \emph{multiplicative} restrictions in~\cite{Yakubovich}. 
There are also restrictions which do not fit into this framework. 
Sometimes, however, there does exist a bijection with a set of partitions which is multiplicative. 

\begin{example}{\rm 
Consider the set of partitions which do not contain any part exactly once, i.e., $R = \{0,2,3,4,\ldots\}$, which is an example of a restriction which is multiplicative in the sense of~\cite{Yakubovich}. 
By conjugation of each element in this set, we obtain the set of partitions which do not contain two consecutive parts and also do not contain any parts equal to 1; see~\cite{RomikUnpublished}. 
We cannot apply Theorem~\ref{integer:partition:theorem} to either of these sets of partitions since the gap at having multiplicity 1 violates the perturbed log-concave condition. 
}\end{example}

\begin{example}{\rm \label{convex:example}
For an integer partition $\lambda$, denote the part sizes in decreasing order as $\lambda_1 \geq \lambda_2 \geq \ldots$, and let $\ell$ denote the number of parts. 
A partition is called \emph{convex} if it satisfies $\lambda_1 -\lambda_2 \geq \lambda_2 - \lambda_3 \geq \ldots \geq \lambda_\ell > 0$; see~\cite{AndrewsII}.
Generalizing this further, for any positive integer $r$, let $\triangle_i^r(\lambda)$ denote the $r$-th order difference, which we define as 
\[ \triangle_i^r(\lambda) =
\begin{cases}
\lambda_i, & \text{if } i=\ell \text{ or } r = 0 \\
\triangle_i^{r-1}(\lambda) - \triangle_{i+1}^{r-1}(\lambda), & \text{otherwise}.
\end{cases}
\]
Consider more generally the set of partitions for which $\triangle_i^r(\lambda) \geq 0$ for all $i$; see~\cite{AndrewsII, canfield2001random}. 
Restrictions of this form a priori fall outside the scope of Theorem~\ref{integer:partition:theorem}; however, in this case we have the following redeeming bijection.
\begin{theorem}[\cite{NonNegR}]\label{convex:theorem}
For all $n=1,2,\ldots$, the number of integer partitions of $n$ which satisfy $\triangle_i^r(\lambda) \geq 0$ is equal to the number of integer partitions of $n$ with unrestricted part sizes in the set $U$ with $u_k = \binom{r+k-1}{r}$, $k=1,2,\ldots$.
\end{theorem}
The mere \emph{existence} of this bijection is what allows us to extend the scope of Theorem~\ref{integer:partition:theorem} to this set of restricted integer partitions; see Proposition~\ref{convex:proposition} and Example~\ref{example:self:conjugate} for examples where the form of the bijection is utilized. 
}\end{example}

\begin{example}{\rm 
Another restriction, known as minimal difference-$d$, imposes that the parts in a partition have difference at least $d$. 
When $d=1$ this is equivalent to partitions into distinct parts, but when $d=2$ it does not fit into the framework of Theorem~\ref{integer:partition:theorem}. 
Nevertheless, there exists a bijection for which we may approximate such partitions asymptotically, which we now describe, see, e.g.,~\cite{Romik, PakSurvey}. 

For any partition $\lambda$ of size~$n$ into exactly $k$ parts, add $d(k-i)$ to $\lambda_i$, $i=1,2,\ldots, k$; call this new partition $\mu$.  
Then $\mu$ is an integer partition of size~$n+d\binom{k}{2}$, into exactly $k$ parts, in which part sizes are at least $d$ apart. 
In other words, we have a bijection between the set of minimal difference-$d$ partitions into exactly $k$ parts and unrestricted integer partitions of $n-d\binom{k}{2}$ into exactly $k$ parts. 
It is also shown in~\cite[Theorem~3]{Romik} that the typical number of parts is asymptotically $\gamma_d \sqrt{n}$, for some $\gamma_d$ which is given in terms of implicitly defined parameters. 
Thus, Theorem~\ref{integer:partition:theorem} can be used to recover~\cite[Theorem~2]{Romik}.  
}\end{example}

There are also more exotic conditions which do not fit into this framework altogether. 
For example, one could demand that each multiplicity is unique, i.e., there do not exist any parts in the partition which occur the same number of times, see~\cite{kane2013asymptotics}. 

\subsection{Partition statistics}
\label{partition:statistics}
The previous section did \emph{not} make use of the Arratia--Pittel--Tavare principle, instead relying only on the quantitative local central limit theorem of Theorem~\ref{p_qlclt_for_combinatorics}. 
We now present a specialization of Theorem~\ref{AT:principle} for integer partitions. 
The special case of unrestricted integer partitions was first achieved in~\cite{PittelShape}.

\begin{theorem}\label{dtv:U:a}
Suppose $U \subset \{1,2, \ldots, n\}$, with $u_k \sim B\, k^r$ for $k \geq 1$, $B >0$, $r \geq 1$.
For any positive integer $t$, suppose $R_i = \{0, 1, \ldots, t\}$ for $i \in U$. 
For any $j_1, j_2$ positive integers, let 
\begin{align}
B_- & := \{1, \ldots, j_1\} \cap U \\
B_{+} & := \{j_2, \ldots, n\} \cap U.
\end{align}
Take $Z_i$ as in Equation~\eqref{eq:Z} with $g_i(k) = \mathbbm{1}(i \in U, k \in R_i)$ and $x = e^{-c_{r,B,t}/n^{r/(r+1)}}$, 
where (see for example~\cite{Ingham, Goh})
\[ c_{r,B,t} = \left\{\frac{(1-(1+t)^{-1/r})}{B^{1/r}\, r} \Gamma\left(\frac{1}{r}+1\right)\zeta\left(\frac{1}{r}+1\right)\right\}^{r/(r+1)}. \]
Then 
\begin{enumerate}
\item $j_1 / n^{r/(r+1)} \to 0$ implies $d_{TV}(Z_{B_-}, C_{B_-}) \to 0$.
\item $j_2 / n^{r/(r+1)} \to \infty$ and $j_2 \leq \chi\,n^{r/(r+1)} \log(n)$, with $\chi < r/((1+r)c_{r,B,t})$, implies $d_{TV}(Z_{B_-}, C_{B_-}) \to 0.$
\end{enumerate}
\end{theorem}

\begin{theorem}\label{dtv:U}
Theorem~\ref{dtv:U:a} continues hold with $R_i = \{0,1,\ldots\}$ by replacing $c_{r,B,t}$ with $\lim_{t\to\infty} c_{r,B,t}$.  
\end{theorem}

The proof of this theorem is a straightforward application of Theorem~\ref{AT:normal}, by an elementary analysis of the random variable $\sum_{i\in U} i\, Z_i$ akin to the first part of Section~\ref{section:asymptotic:enumeration}. 
In particular, this theorem tells us that the set of integer partitions into distinct parts follows a remarkable similarity to unrestricted integer partitions, which was also observed in~\cite{Fristedt}.  

Theorem~\ref{dtv:U} with $r=B=1$ was utilized in~\cite{PittelShape} to prove the limit shape for integer partitions, which requires more analysis on the ``middle" region with ``moderately"--sized parts. 
At present, the theory of limit shapes is sufficiently mature in that other related approaches have been utilized to obtain such a statistic in many general settings; see for example~\cite{Yakubovich, LD, Goh, Bogachev}. 
Thus, since one does not immediately obtain the limit shape from our theorems, i.e., it requires a further detailed analysis of the middle region, where total variation distance does not tend to zero, we do not pursue this further. 

There are, however, a plethora of interesting statistics on integer partitions which are localized on the smallest and largest part sizes, and hence are well-approximated asymptotically using Theorem~\ref{dtv:U:a}. 
In the examples that follow, we have chosen to assume the restrictions in Theorem~\ref{dtv:U}, since it is the type that comes up most often in our applications and very nicely generalizes many classical results. 

\begin{example}{\rm 
A \emph{gap} in an integer partition, see also Example~\ref{example:gap}, is a part size which does not occur in the partition; i.e., a gap occurs at $i$ when $C_i(n) = 0$. 
With the assumptions in Theorem~\ref{dtv:U}, the probability of a gap at index~$u_1$ is given by 
\[ \P(C_{u_1}(n) = 0) \sim \P(Z_{u_1} = 0) = 1-x^{u_1} = 1-e^{-u_1\, c_{r,B}/n^{r/(r+1)}} \sim \frac{u_1\, c_{r,B}}{n^{r/(r+1)}}. \]
Recall that a special case of this was utilized in Equation~\eqref{eq:no:one} to approximate the number of partitions not containing any parts of size~$1$. 
Similarly, for any $i\in U$ such that $i = o(n^{r/(r+1)})$,  we have
\[ \P(C_i(n) = 0) \sim \P(Z_i = 0) = 1-x^i = 1-e^{-c_{r,B}\,i /n^{r/(r+1)}} \sim \frac{c_{r,B}\, i}{n^{r/(r+1)}}. \]
One could ask, then, what is the typical size of the smallest gap in a random integer partition of size~$n$; see~\cite{grabner2006analysis}. 
Let $I := \min(i : Z_i = 0)$.  For $k = o(n^{r/(r+1)})$, we have 
\[ \P(I = k) \sim \P(Z_k = 0) \prod_{i=1}^{k-1}\P(Z_i \geq 1)  = (1-x^k) \prod_{i=1}^{k-1} x^i = (1-x^k) x^{\binom{k}{2}}. \]
From the term $x^{\binom{k}{2}}$, we see that the majority of the contribution to the random variable $I$ is from values of $k = O(\sqrt{n^{r/(r+1)}})$. 
Thus, the expected value is then asymptotically 
\[ \e I \sim \sum_{k=1}^{\sqrt{n^{r/(r+1)+\epsilon}}} k (1-x^k)x^{\binom{k}{2}}\mathbbm{1}(k \in U) \sim \sum_{k\geq 1} u_k \left(1-\exp\left({-c_{r,B} u_k/n^{r/(r+1)}}\right)\right) \exp\left({-c_{r,B}\frac{\binom{u_k}{2}}{n^{r/(r+1)}}}\right). \]
Define $y_k$ such that $B\, y_k^r = u_k / n^{r/(2(r+1)}$, so that $\Delta y_k = y_k - y_{k-1} \sim 1/n^{1/(2(r+1))}$. 
We then have 
\begin{align}
 \e I & \sim \sum_{y_k \geq 0} B\, y_k^r\, n^{r/2(r+1)} \left(1-\exp\left(-c_{r,B} B\, y_k^r / n^{r/2(r+1)}\right)\right) \exp\left({-c_{r,B}\frac{y_k^2/2}{n^{r/(r+1)}}}\right) \Delta y_k\ n^{1/2(r+1)} \\
    & \sim  n^{(r+1)/2(r+1)} \int_0^\infty B y^r \left(\frac{c_{r,B} B y^r}{n^{r/2(r+1)}}\right) \exp\left(-c_{r,B} \frac{B^2 y^{2r}}{2}\right) \, dy \\
    & \sim n^{\frac{1}{2} - \frac{r}{2(r+1)}} \int_0^\infty c_{r,B}\, B^2 y^{2r} \exp\left(-\frac{c_{r,B} B^2}{2} y^{2r}\right) dy \\
    & \sim n^{\frac{1}{2(r+1)}} \left(\frac{2}{B\, c_{r,B}}\right)^{\frac{1}{2r}}\ \frac{\Gamma\left(1+\frac{1}{2r}\right)}{r} \ .
\end{align}
In the special case $r=1$ and $B=1$, we recover the latter part of~\cite[Theorem~3]{goh1992gap}, i.e., that the average size of the smallest gap is asymptotically $\sqrt[4]{3n/2}$. 
This value was obtained via a careful analysis of the relevant generating functions. 
Using our approach, we can easily obtain arbitrary moments via a similar calculation, viz.,
\begin{align}
 \e I^s & \sim n^{\frac{r\,s+1-r}{2(r+1)}} \int_0^\infty c_{r,B}\, B^{s+1}\, y^{r\,s}  \exp\left(-\frac{c_{r,B} B^2}{2} y^{2r}\right) dy \\
 &  \sim n^{\frac{r\,s+1-r}{2(r+1)}} c_{r,B}\, B^{s+1}\, \frac{2^{\frac{1}{2} \left(\frac{1}{r}+s-2\right)} B^{-\frac{r s+1}{r}} c_{r,B}^{-\frac{r s+1}{2 r}} \Gamma \left(\frac{r s+1}{2 r}\right)}{r} , \qquad s \geq 1. 
\end{align}
}\end{example}

\begin{remark}
Aside from the small part sizes, one may also inquire as to the behavior of the large part sizes, as was investigated in~\cite{Fristedt}, where it was shown that for unrestricted integer partitions the largest part sizes can be described asymptotically as a Markov chain. 
We have not carried out the detailed analysis, however, it is highly plausible that under the assumptions of Theorem~\ref{dtv:U}, that the largest part sizes also follow an analogous Markov chain law, from which one might provide a generalization to many of the results in~\cite{PittelConfirming}.
\end{remark}

It is also of importance to note that the total variation metric is robust, in the sense that if $h$ is any measurable function, we have 
\begin{align}
\label{dtv:h} d_{TV}( h(C_B), h(Z_B) ) \leq d_{TV}( C_B, Z_B); &
\end{align}
i.e., applying any measurable function to the joint distributions cannot increase the total variation distance. 

\begin{example}{\rm 
Suppose $u_k = \binom{k+1}{2}$, $k=1,2,\ldots$; i.e., partitions into ``triangular" parts. 
Recall the definition of a convex partition, which is defined in terms of its \emph{part sizes} $\lambda_1 \geq \lambda_2 \geq \ldots \lambda_\ell>0$, which satisfy the condition $\lambda_1 - \lambda_2 \geq \lambda_2 - \lambda_3 \geq \ldots \geq \lambda_\ell > 0$.  
Theorem~\ref{convex:theorem} with $r=2$ implies that the number of integer partitions of $n$ into triangular parts is equal to the number of integer partitions of $n$ which are convex, see~\cite{AndrewsOmega}. 
This theorem was of particular interest in Example~\ref{convex:example}, where it extended the scope of Theorem~\ref{integer:partition:theorem} via the \emph{existence} of a bijection. 
However, it is the form of the bijection itself which is most enticing for us in this section. 
That is, letting $a_k$ denote the number of parts of size $\binom{k+1}{2}$, $k=1,2,\ldots$, in a partition into triangular parts, the bijection in~\cite{AndrewsOmega} (see also~\cite{NonNegR, PakSurvey}) creates a convex partition via
\begin{align}
 \lambda_1 & = \sum_{j \geq 1} j\, a_j, \qquad  \lambda_2  = \sum_{j \geq 2} (j-1)\, a_j,  \qquad \cdots \qquad  \lambda_k  = \sum_{j \geq k} (j-k+1)\, a_j. 
\end{align}
In fact, this bijection is simply a linear transformation of multiplicities, i.e., a \emph{measurable function}, and as such we conclude that the joint distribution of smallest part sizes in a convex partition are approximately a weighted sum of certain geometric random variables. 
A more precise statement is below. 
\begin{proposition}\label{convex:proposition}{\rm 
Let $\lambda_k$ denote the $k$-th largest part size in a random convex partition, and let $X_k := \sum_{j \geq k} (j-k+1)\, Z_{\binom{j+1}{2}}$, $k=1,2,\ldots$. 
Suppose $k_n = \omega(n^{2/3})$ and $k_n \leq \frac{2}{3\, c_{1,1/2}}\ n^{2/3} \log(n)$. 
We have 
\[ d_{TV}\left( (\lambda_{k_n}, \lambda_{k_n+1}, \ldots), (X_{k_n}, X_{k_n+1},\ldots)\right) \to 0. \]
}\end{proposition} 
}\end{example}

\begin{example}{\rm \label{dtv:root:two}
It is, of course, not necessary that we consider contiguous sets of indices.  
Consider the set $B = \{1, 3, 5, \ldots \}$, i.e., the set of all odd part sizes in an integer partition of size~$n$. 
We have 
\[ d_{TV}( C_B, Z_B) \to d_{TV}(\mathcal{N}_B, \mathcal{N}), \]
where $\mathcal{N}_B$ is a normal distribution with variance given by $1/2\sqrt{2}$.  
Thus, if we wish to approximate the odd part-sizes, we have a non-trivial, but also asymptotically non-zero, total variation distance error. 
This was exploited in~\cite{PDC} for an asymptotically efficient random sampling algorithm of integer partitions. 
}\end{example}

\section{Assemblies, multisets, selections}
\label{section:classes}

\subsection{Overview}

Integer partitions are a specialization of a broader class of combinatorial structures. 
In~\cite{IPARCS}, three main classes of combinatorial structures: assemblies, multisets, and selections, are presented in a unified probabilistic manner via specializations of Equation~\eqref{eq:Z}.

A generalization to classical integer partitions which has combinatorial meaning occurs when taking $g_i(k) = \binom{m_i+k-1}{k}$ in Equation~\eqref{eq:Z}, for ${\bf m} \equiv (m_1, m_2, \ldots)$ nonnegative integers, from which we have 
\begin{align}
\label{mi:partitions} \P\left(\sum_{i=1}^n i\, Z_i = n \right)  & = \sum_{(c_1, \ldots, c_n) \in \mathbb{Z}_0^n : \sum_{k} k\, c_k = n} \binom{m_k-c_k-1}{c_k} x^{n} \prod_{k=1}^n(1-x^k)^{m_k} \\
&  = p_{\bf m}(n)\, x^{n} \prod_{k=1}^n \left(1-x^k\right)^{m_k}, 
\end{align}
where $p_{\bf m}(n) := \sum_{(c_1, \ldots, c_n) \in \mathbb{Z}_0^n : \sum_{k} k\, c_k = n} \binom{m_k-c_k-1}{c_k}$. 
Note that, unlike in the case of integer partitions, there is a weight attached to each summand.
This weight assigns the appropriate proportion among all objects with that given component structure. 
The random variables $Z_k$ have a \mbox{NegativeBinomial} distribution with parameters $m_k$ and $x^k$, $k=1,2,\ldots,n$, and the quantity $p_{\bf m}(n)$ can now be approximated assuming the left-hand side~of~\eqref{mi:partitions} satisfies a LCLT. 
This particular family of distributions corresponds to a set of combinatorial structures known as \emph{multisets}. 

Another important family of distributions arises when $g_k(\ell) = \binom{m_k}{\ell}$ in Equation~\eqref{eq:Z}, $k=1,2,\ldots,$ $\ell=0,1,\ldots$, which correspond to \emph{selections}. 
In this case, we have $c_k = (1-x^k)^{-m_k}$, and hence the $Z_k$ have a \mbox{Binomial} distribution with parameters $m_k$ and $\frac{x^k}{1+x^k}$. 
By writing the binomial distribution as a sum of independent Bernoulli random variables, we can think of $\sum_{k=1}^n k\, Z_k = \sum_{k=1}^n k\, (B_{1,k} + \cdots + B_{m_k, k})$ as selecting a component of weight $k$ independently $m_k$ times, each with probability $x^k/(1+x^k)$, $k=1,2,\ldots, n$, and labeling the chosen components according to which Bernoulli random variable generated them. 
Letting $s_{\bf m}(n)$ denote the number of selections, we have 
\begin{align}
\label{mi:selections} \P\left(\sum_{i=1}^n i\, Z_i = n \right) = s_{\bf m}(n)\, x^{n} \prod_{k=1}^n \left(1+x^k\right)^{-m_k}. 
\end{align}

Our final notable family of distributions arises when $g_k(\ell) = (m_k/k!)^\ell/\ell!$ in Equation~\eqref{eq:Z}, corresponding to \emph{assemblies}, implying $Z_k$ is Poisson with parameter $\lambda_k = m_k\ts x^k / k!$, $k=1,2,\ldots, n$. 
This family of distributions corresponds to the blocks of a random set partition when $m_k = 1$ for all $k=1,2,\ldots,n$, and to the cycle lengths of a random permutation when $m_k = (k-1)!$. 
Letting $a_{\bf m}(n)$ denote the number of assemblies, we have 
\begin{align}
\label{mi:assemblies} \P\left(\sum_{i=1}^n i\, Z_i = n \right) = a_{\bf m}(n)\, \frac{x^{n}}{n!} \exp\left(- \sum_{k=1}^n \frac{m_k x^k}{k!}\right). 
\end{align}

\begin{remark}
It is also possible to generalize this probabilistic framework further, say to linear combinations of the random variables, i.e., $\sum_{k=1}^n w_k Z_k$, see for example~\cite[Section~2.2]{IPARCS}, but we do not pursue this presently. 
\end{remark}

\subsection{Special cases with known LCLTs}

Before proceeding with our extension of theorems~\ref{p_qlclt_for_combinatorics} and~\ref{AT:principle}, we first review the existing results governing sums of the form $\sum_{k=1}^n k\, Z_k$ for the families of distributions satisfying~\eqref{eq:Z}. 

The first large class refers to \emph{logarithmic combinatorial structures}, see~\cite{arratia1995total, Logarithmic}, which only require the following two conditions on the random variables $Z_k$ composing the sum $\sum_{k=1}^n k\ts Z_k$: for some $\theta >0$, as $k \to \infty$ we have 
\begin{align}
& k\, \P(Z_k = 1) \to \theta, \\
& k\, \e Z_k \to \theta. 
\end{align}
As pointed out in~\cite[Section~3.3]{Logarithmic}, see also~\cite{IPARCS}, there are corresponding growth rates of $m_k$ in each of the three classes of families for which the above logarithmic condition holds for fixed $x$: 
\[\begin{array}{lll}
\mbox{assemblies:} & m_k \sim \kappa\ts y^k (k-1)!, & \kappa > 0, y>0, x = y^{-1}; \\
\mbox{multisets:} & m_k \sim \frac{\kappa\ts y^k}{k}, & \kappa > 0, y>1, x = y^{-1}; \\
\mbox{selections:} & m_k \sim \frac{\kappa\ts (1+y^k)}{k}, & \kappa > 0, y>1, x = y^{-1}.
\end{array} \]
Under these assumptions, it was established in~\cite[Chapter~5]{Logarithmic} that indeed these three classes satisfy a LCLT \emph{which is not necessarily Gaussian}. 

A less restrictive set of conditions for the LCLT to hold in the case $m_k = 1$ for all $k$ are those contained in~\cite{Bogachev}, motivated by the finding of limit shapes of these combinatorial families. 
The LCLT is the most technical part of establishing such a limit shape, and even then one need only find an appropriately sub-exponential lower bound, as is the approach in, e.g.,~\cite{RomikUnpublished}. 
As stated in the introduction, the conditions imposed for a LCLT to hold in~\cite{Bogachev} require an analysis of the generating function. 
Specifically, one takes the logarithm of the generating function and writes out the Taylor series expansion, and assumes another technical condition which we now specify. 
Denote by $F(x)$ the generating function and $H(x) := \log(F(x)) = \sum_{k \geq 0} a_k x^k$.  Then if there exists a constant $\delta_\ast$ such that for any $\theta \in (0,1)$ we have 
\begin{equation}\label{Bogachev:smooth} H(\theta) - Re( H(\theta e^{i\ts t})) \geq \delta_\ast \theta (1-\cos t), \qquad t \in \mathbb{R}, \end{equation}
and 
\[ \sum_{k \geq 0} \frac{|a_0|}{\sqrt{k}} < \infty, \]
then a LCLT holds for the sum $\sum_{k=1}^n k\, Z_k$. 
Inequality~\eqref{Bogachev:smooth} is an application of $H$-admissibility, see~\cite[Section~12.2]{odlyzko1995asymptotic}, which we have replaced with our condition of perturbed logconcavity on the distributions in~\eqref{eq:Z}, and the finding of the set $M$ corresponding to the peak in variances. 
Conditions like~\eqref{Bogachev:smooth} are quite common in the theory of integer partitions, see for example~\cite[Condition~(II)]{RothSzekeres} and~\cite[Theorem~1]{Goh} and~\cite{Ingham}. 
To reiterate, even when there is overlap in the existing theory, our conditions are often easier to verify in practice. 

Another general treatment for limit shapes is~\cite{Yakubovich}, which also contains a local central limit theorem.  
The conditions for~\cite[Lemma~10]{Yakubovich} to hold involve the radius of convergence of the generating function of the sequence, a periodicity condition akin to Equation~\eqref{Bogachev:smooth}, and a control on the rate of increase of $m_k$. 

Another recent work is~\cite{erlihson2008limit}, which can be specialized and applied \emph{for assemblies only}, assumes for any $C>0$ and $p \geq 0$, we take $m_k \sim C k^{p-1} k!$. 
This restriction places these structures in the realm of logarithmic combinatorial structures, see~\cite{Logarithmic}. 
In this setting, a local central limit theorem is obtained as well as a bivariate local central limit theorem for the event $\{\sum_i i\, Z_i, \sum_i Z_i\}$. 
This is used to also obtain a joint distribution on part sizes similar to that in~\cite{sachkov1974random}, i.e., a jointly \emph{normal} distribution. 
Similar to~\cite{Bogachev}, the motivation is to obtain a limit shape, for which the local central limit theorem is a technical tool which, combined with a large deviation result, implies the equivalence of the limit shapes of the dependent and independent processes (i.e., the micro and macro canonical ensembles); see also~\cite{Yakubovich}. 

\subsection{Multisets}
A multiset is a generalization to an integer partition.  
The parameters $m_i$ determine the number of different types of components of size~$i$, for $i=1,2,\ldots, n$. 
A powerful theorem of~\cite{Meinardus} is applicable in this setting, which covers a wide range of possible values for $m_i$, though there are examples where this theorem is not directly applicable, see~\cite{romik2015number}. 
It requires several definitions which we simply state, and refer the interested reader to~\cite[VIII.23]{Flajolet} or~\cite[Chapter~6]{Andrews}.

\begin{theorem}[\cite{Meinardus}, see also~\cite{Flajolet, Andrews}]\label{meinardus:theorem}
Suppose $a_1, a_2, \ldots$ are given nonnegative (not necessarily integer) values.
Define $f(z)$ and $f_1, f_2, \ldots$ as 
\[ f(z) := \prod_{i=1}^\infty (1-z^i)^{-a_i} = \sum_{i \geq 1} f_i z^i. \] 
Let $\alpha(s) := \sum_{i\geq 1} \frac{a_i}{i^s}$ denote the associated Dirichlet series, and assume that $\alpha(s)$ is continuable into a meromorphic function to $\mathcal{R}(s) \geq -C_0$ for some $C_0$, with only a simple pole at some $\rho>0$ and corresponding residue $A$. 
Assume also for some $C_1>0$, as $|s| \to \infty$ such that $\mathcal{R}(s) \geq -C_0$, we have $\alpha(s) = O(|s|^{C_1})$. 
Let $g(z) := \sum_{i \geq 1} a_i z^i$, and assume
\[ \mathcal{R} g(e^{-t-2i\pi y}) - g(e^{-t}) \leq - C_2 y^{-\epsilon}. \]
Then we have 
\[ f_n = C n^{\kappa} \exp\left( n^{\alpha/(\alpha+1)} \left(1+\frac{1}{\alpha}\right)\left(A\Gamma(\alpha+1)\zeta(\alpha+1)\right)^{1/(1+\alpha)}\right)\left(1+O(n^{-\kappa_1})\right), \]
where
\[ C = e^{\alpha'(0)}\left(2\pi(1+\alpha)\right)^{-\frac{1}{2}} \left(A \Gamma(\alpha+1)\zeta(\alpha+1)\right)^{(1-2\alpha(0))/(2+2\alpha)}, \]
\[ \kappa = \frac{\alpha(0) - 1 - \frac{1}{2}\alpha}{1+\alpha}, \]
\[ \kappa_1 = \frac{\alpha}{\alpha+1} \min\left(\frac{C_0}{\alpha}-\frac{\delta}{4}, \frac{1}{2}-\delta\right), \]
where $\delta$ is an arbitrary real number. 
\end{theorem}

Despite the power of the previous theorem, it is our sincere desire to forego the standard analytical assumptions on generating functions and instead embrace a more simplistic probabilistic formulation. 
We next state a fairly broad specialization of Theorem~\ref{AT:principle} to the parameterization $m_k = a_d k^d + a_{d-1}k^{d-1}+ \ldots+ a_0$, along with the restriction that the component sizes come from the set $U$ satisfying the assumptions in Theorem~\ref{ingham:theorem}. 

\begin{theorem}\label{multisets:theorem}
Suppose $U = \{u_1, u_2, \ldots\}$ is a given subset of positive integers which satisfies the assumptions in Theorem~\ref{ingham:theorem}, namely, that $N(u)$ denotes the number of elements in $U$ less than or equal to $u$, and there are constants $A>0$, $\theta >0$, and function $R(u)$ such that 
\[ N(u) = A\, u^\theta + R(u), \]
and for some $b>0$ and $c>0$ we also have 
\[ \int_0^u \frac{R(v)}{v}\, dv = b \log u + c + O(f(u)), \]
where $f(u) = o(1)$. 
Suppose ${\bf m} = (m_1, m_2, \ldots )$ is a given collection of nonnegative integers such that $m_k = a_d k^d + a_{d-1}k^{d-1}+ \ldots+ a_0$, $a_d \neq 0$. 
Let $p(n; U,{\bf m})$ denote the number of multisets of $n$ with component sizes in $U$, with $m_k$ components of size~$k$. 
Then we have 
\begin{align}
\label{multisets:asymptotic}p(n; U, {\bf m}) & = \frac{e^{-\alpha\, n} \prod_{i\in U, i \leq n}\left(1-e^{-\alpha\, i}\right)^{-m_k}}{\sqrt{2\pi\sigma^2}}\left(1+O\left(n^{-\frac{d+\frac{1}{r}}{2(d+1+\frac{1}{r})}}\right)\right),
\end{align}
where $r := 1/\theta$, 
\[ \alpha = \left(\frac{c_1}{n}\right)^{\frac{r}{r\,d+r+1}}, \qquad \qquad \sigma^2 = n^{\frac{d+2+1/r}{d+1+1/r}} \, c_2, \]
\[ c_1 = \int_0^\infty a_d B^{d+1} y^{r(d+1)} \frac{e^{-B\, y^r}}{1-e^{-B\,y^r}}\, dy, \]
\[ c_2 =  \int_0^\infty a_d\, B^{d+2}\, y^{r(d+2)}\, \frac{e^{-By^r}}{(1-e^{-By^r})^2}\, dy. \] 
\end{theorem}

Of course, Theorem~\ref{multisets:theorem} still requires careful analysis of the product $\prod_{i\in U} \left(1-e^{-\alpha\, i}\right)^{-m_k}$, see for example the proofs of~\cite[Proposition~1]{Romik},~\cite[Theorem~2]{Ingham}. 
For example, we recover~\cite[Lemma~6]{Dalal} using $U = \{1,2,\ldots\}$ and $m_k=2$ for all $k$. 
For the approximation theorem below, we shall not need the full asymptotic formula. 

When $m_k = k$ for all $k=1,2,\ldots,$ we may interpret this structure combinatorially as the set of integer partitions where each part of size~$k$ can be colored in any of $k$ distinct colors. 
It is also the generating function for the number of plane partitions of $n$, as pointed out in~\cite{MacMahon}. 

\begin{theorem}[\cite{WrightPlane}, see also~\cite{Flajolet, BodiniPlane, Andrews}]\label{plane:partitions:theorem}
Let $PP(n)$ denote the number of plane partitions of size~$n$. 
We have 
\[ PP(n)  = \frac{(\zeta(3)\, 2^{-11})^{1/36}}{n^{25/36}} \exp\left(3\, 2^{-2/3} \zeta(3)^{1/3} n^{2/3} + 2d\right)\left(1+O(n^{-1/3})\right), \]
where $d = -\frac{e}{4\pi^2} \left(\log(2\pi) + \gamma - 1\right).$
\end{theorem}

An equivalent approach for proving Theorem~\ref{plane:partitions:theorem} above is in~\cite{BodiniPlane}, which uses an array of geometric random variables as follows. 
For $j,k \geq 0$, let $Z_{j,k}$ denote geometric random variables with parameter $1-x^{j+k+1}$, and consider
$S_n = \sum_{i,j\geq 0} (j+k+1) Z_{j,k}.$ 
Then we have
\[ \p(S_n = n) = f(n) x^n \prod_{j,k \geq 0} (1-x^{j+k+1}). \]
A priori, it is not apparent that $f(n)$ above has any classical combinatorial meaning; however, the bijection given in~\cite{PakBijection} shows that indeed $f(n) = PP(n)$. 
We then perform the same calculations as in previous examples. 
Solving for $\e S_n \sim n$, we let $x = e^{- (2\zeta(3) / n)^{1/3}}$, where $\zeta$ denotes the Riemann zeta function, $\zeta(r) = \sum_{k\geq 1} r^{-k}$.  
We then have
\[ PP(n) = \frac{e^{n^{2/3}(2\zeta(3))^{1/3}} \prod_{j,k \geq 0} (1-e^{- (j+k+1)\, (2\zeta(3) / n)^{1/3}})^{-1}}{ \sqrt{2\pi} \sigma}\left(1 + O(n^{-1/3}) \right). \]
This approach is advantageous in that it readily generalizes to plane partitions conditioned to lie in some box $[0,0] \times [b_1,b_2]$ in the plane, see~\cite{BodiniPlane}. 
In this case, the sum of interest is 
\[ \sum_{j=1}^{b_1} \sum_{k=1}^{b_2} (j+k+1) Z_{j,k}, \]
and the rest of the application mimics the calculations above. 

Theorem~\ref{AT:principle} is specialized as follows. 

\begin{theorem}\label{multiset:dtv}
Suppose $U = \{u_1, u_2, \ldots\}$ is a given subset of positive integers such that $u_k \sim b\, k^r$. 
Suppose ${\bf m} = (m_1, m_2, \ldots )$ is a given collection of nonnegative integers such that $m_k \sim a\, k^\theta$. 
For any $j_1, j_2$ positive integers, let 
\begin{align}
B_- & := \{1, \ldots, j_1\} \cap U \\
B_{+} & := \{j_2, \ldots, n\} \cap U.
\end{align}
Take $Z_i$ as in Equation~\eqref{eq:Z} with $g_i(k) = \binom{m_i+k-1}{k}\, \mathbbm{1}(i \in U)$ and $x = e^{-c_{1}/n^{r/(rd+r+1)}}$, where 
\[ c_1 = \left(\int_0^\infty a_d B^{d+1} y^{r(d+1)} \frac{e^{-B\, y^r}}{1-e^{-B\,y^r}}\, dy\right)^{r/(rd+r+1)}. \]
Then 
\begin{enumerate}
\item $j_1 / n^{r/(rd+r+1)} \to 0$ implies $d_{TV}(Z_{B_-}, C_{B_-}) \to 0$.
\item $j_2 / n^{r/(rd+r+1)} \to \infty$ and $j_2 \leq \chi\,n^{r/(rd+r+1)} \log(n)$, with $\chi < r/((rd+r+1)c_{1})$, implies $d_{TV}(Z_{B_-}, C_{B_-}) \to 0.$
\end{enumerate}
\end{theorem}

\begin{remark}{\rm
In contrast to, e.g., Example~\ref{example:self:conjugate}, it is not known whether there exists a bijection for plane partitions which corresponds to a continuous transformation of the component sizes $\{Z_{i,j}\}_{1 \leq i, 1 \leq j}$, in the sense described in~\cite{Pak2}. 
If such a continuous bijection were discovered, one could then apply Theorem~\ref{multiset:dtv} via Equation~\eqref{dtv:h} in order to describe various statistics on plane partitions. 
}\end{remark}

\subsection{Selections} 
\label{selections}
The class of selections are those combinatorial objects for which components of size~$k$ are either selected or not selected, hence their values lie in the set $\{0,1\}$. 
Selections follow the same form as multisets, in that there are parameters $m_1, m_2, \ldots$ which we say define the selection, and $m_k$ denotes the number distinguishable components of size~$k$, and each distinguishable component can only be represented at most once in the object. 
Letting $X_k$ have binomial distribution with parameters $m_k$ and $\frac{x^k}{1+x^k}$, for any $0<x<1$, and for any combinatorial selection with $c_k$ parts of size~$k$, $k=1,2,\ldots,$ with $\sum_k k\, c_k = n$, we have 
\[ \p(X_1 = c_1, \ldots, X_n = c_n) = \prod_{k=1}^n \binom{m_k}{c_k} \left(\frac{x^k}{1+x^k}\right)^{c_k}\left(\frac{1}{1+x^k}\right)^{m_k-c_k}.\]
Summing over all possible selections of weight~$N$, and simplifying using the Binomial theorem, we have
\[ \p\left(\sum_{k=1}^n k\, X_k = n\right) = s(n) x^n \prod_{k=1}^n (1+x^k)^{m_k},\]
where $s(n)$ is the total number of selections of weight~$n$.

\begin{theorem}\label{selections:theorem}
Suppose $U$ and ${\bf m}$ satisfy the conditions of Theorem~\ref{multisets:theorem}. 
Take $Z_i$ as in Equation~\eqref{eq:Z} with $g_i(k) = \binom{m_i}{k}\, \mathbbm{1}(i \in U)$ and $x = e^{-d_{1}/n^{r/(rd+r+1)}}$, 
where 
\[ d_1 = \left(\int_0^\infty a_d B^{d+1} y^{r(d+1)} \frac{e^{-B\, y^r}}{1+e^{-B\,y^r}}\, dy\right)^{r/(rd+r+1)}. \]
Let $s(n; U,{\bf m})$ denote the number of selections of $n$ with component sizes in $U$, with $m_k$ components of size~$k$. 
We have 
\[ s(n; U, {\bf m}) = \frac{e^{-\alpha\, n} \prod_{i\in U, i \leq n}\left(1+e^{-\alpha\, i}\right)^{m_k}}{\sqrt{2\pi\sigma^2}}\left(1+O\left(n^{-\frac{d+\frac{1}{r}}{2(d+1+\frac{1}{r})}}\right)\right), \]
where $r := 1/\theta$, 
\[ \alpha = \left(\frac{d_1}{n}\right)^{\frac{r}{r\,d+r+1}}, \qquad \qquad \sigma^2 = n^{\frac{d+2+1/r}{d+1+1/r}} \, d_2, \]
\[ d_1 = \int_0^\infty a_d B^{d+1} y^{r(d+1)} \frac{e^{-B\, y^r}}{1+e^{-B\,y^r}}\, dy, \]
\[ d_2 =  \int_0^\infty a_d\, B^{d+2}\, y^{r(d+2)}\, \frac{e^{-By^r}}{(1+e^{-By^r})^2}\, dy, \] 
\end{theorem}

For example, taking $m_k =1$ corresponds to integer partitions into \emph{distinct part sizes}, which follows a similar analysis as in Section~\ref{section:integer:partitions}. 
We also note that the same analysis for the product $g(\alpha) := \prod_{i\in U} \left(1-e^{-\alpha\, i}\right)^{-m_k}$ in Theorem~\ref{multisets:theorem} can be utilized in Theorem~\ref{selections:theorem} directly, from the fact that $\prod_{i \in U, i\leq n}\left(1+e^{-\alpha\, i}\right)^{m_k} = g(2\alpha)/g(\alpha)$. 

\begin{example}{\rm \label{example:self:conjugate}
An example which does \emph{not} immediately follow from previous work is the following. 
For an integer partition $\lambda = (\lambda_1, \ldots, \lambda_\ell)$ of size~$n$, with part sizes $\lambda_1 \geq \lambda_2 \geq \cdots \geq \lambda_\ell > 0$, the conjugate partition $\lambda'$ is defined by $\lambda_j' = \sum_{i \geq j} \mathbbm{1}(\lambda_i \geq j)$. 
For example, $\lambda_1' = \ell$, the total number of parts in $\lambda$. 
A partition is called \emph{self-conjugate} whenever $\lambda = \lambda'$. 

The limit shape of self-conjugate partitions is known to coincide with the limit shape of unrestricted integer partitions, even though the number of such integer partitions is of a different exponential asymptotic order. 
Several combinatorial proofs are contained in~\cite[Section~7]{Pak2}, which are defined in terms of continuous transformations, from which it follows as well that the number of self-conjugate partitions is equal to the number of integer partitions into distinct, odd part sizes. 
Of the bijections presented, the one via hook lengths is the most appropriate for our purpose. 
That is, let 
\begin{equation}\label{hook:bijection} 
\mu_i = 2(\lambda_i - i)+1, \qquad \mbox{for $1 \leq i \leq \delta_0(\lambda)$},\end{equation}
where $\delta_0(\lambda)$ is the largest Durfee square of $\lambda$, which is defined as the largest square that can fit inside the Young diagram in $\lambda$ (see~\cite{Pak2} for the relevant definitions, which we shall not need). 
It can be shown that for partitions into distinct, odd parts, we have $\delta_0(\lambda) \sim (\ln(2)/c) \sqrt{n}$. 

Letting $\lambda' = \mu$ for each $\lambda$ into distinct, odd parts, we see that Equation~\eqref{hook:bijection} defines a continuous transformation of part sizes. 
Thus, by Equation~\eqref{dtv:h}, the same theorems governing partitions into distinct, odd parts, i.e., Theorem~\ref{integer:partition:theorem} and Theorem~\ref{dtv:U} using $U = \{1, 3, 5, \ldots\}$, $R = \{0,1\}$, are applicable to self-conjugate partitions via Equation~\eqref{hook:bijection}. 
}\end{example}

\subsection{Assemblies} 
\label{assemblies}
Our final combinatorial class is known as assemblies, and represents a class of combinatorial structures whose component sizes are well-approximated by Poisson random variables. 
As in the case of multisets and selections, we define a particular assembly using $m_1, m_2, \ldots$. 
In this case, for any $x>0$, the random variable $Z_i$ in Equation~\eqref{eq:Z} is a Poisson random variable with parameter $\lambda_i := m_i x^i / i!$, $i=1,2,\ldots$. 

Taking $m_k = 1$ for all $k=1,2,\ldots,$ the corresponding combinatorial structure is the block sizes of a random set partition, and we have $a_{\bf m}(n)$ is the $n$-th Bell number; see for example~\cite{PittelSetPartitions}.
Taking $m_k = (k-1)!$ for all $k=1,2,\ldots,$ the corresponding combinatorial structure is the cycle sizes of a random permutation of $n$, and we have $a_{\bf m}(n) = n!$. 

Unfortunately, Theorem~\ref{p_qlclt_for_combinatorics} is not strong enough for set partitions, which we now demonstrate. 
A set partition of size~$n$ is a partition of the set $\{1, 2, \ldots, n\}$. 
The subsets are called blocks, and similar to integer partitions, one can describe the block structure of a set partition in terms of the sizes of each block. 
\emph{Unlike} integer partitions, however, is that once we have specified the block sizes, there are many ways to fill them in with numbers from the set $\{1, \ldots, n\}$. 
Thus, while set partitions share many similar properties as integer partitions, it is not surprising that there are key differences.

Taking $g_i(k) = \frac{1}{(i!)^k\, k!}$ in Equation~\eqref{eq:Z}, for any $x>0$ we have 
\[ \P(Z_i = k) = c_i \frac{(x^i/i!)^k}{k!} = e^{-x^i/i!} \frac{(x^i/i!)^k}{k!},  \qquad k=0,1,\ldots. \]
Thus, $Z_i$ has a Poisson distribution with parameter $\lambda_i := x^i/i!$, $i=1,2,\ldots$. 
Letting $T := \sum_{i=1}^n i\, Z_i$ denote the sum of all block sizes, we have 
\[ \P(T = n) = B(n) \frac{x^n}{n!} \exp\left(-\sum_{i=1}^n \frac{x^i}{i!}\right), \]
where $B(n)$ denotes the $n$-th Bell number, which satisfies $B(0) = 1$ and the recurrence
\[ B(n+1) = \sum_{k=0}^n \binom{n}{k} B(k), \qquad n \geq 1. \]

The expected value of $T$ is given by 
\[ \e T = \sum_{i=1}^n i\, \frac{x^i}{i!} = x \sum_{i=0}^{n-1} \frac{x^i}{i!} \sim x e^x. \]
Thus, we wish to take $x>0$ such that $x e^x = n$. 
This solution is known as the product logarithm, or the Lambert W function, and it is typically denoted by $W(n)$.  
In~\cite{Debruijn}, it is shown that 
\[ W(n) = \log(n) - \log\log(n) + \frac{\log(\log(n))}{\log(n)}+ \frac{1}{2}\left(\frac{\log(\log(n))}{\log(n)}\right)^2 + O\left(\frac{\log\log(n)}{\log^2(n)}\right), \]
and we shall often take the first term, i..e, $x \sim \log(n)$ in what follows, although in this case we must exercise more care since the expressions involving $x$ often appear in exponents. 
With $x$ chosen in this manner, the variance of $T$ is given by 
\[ \var{T} = \sum_{i=1}^n i^2 \frac{x^i}{i!} = \sum_{i=1}^n i(i-1) \frac{x^i}{i!} + i \frac{x^i}{i!} \sim x^2 e^x + x e^x \sim n(1+x) \sim n \log(n). \]
More precisely, we have 
\[ |\e T - n| \leq \sum_{k=n+1}^\infty \frac{x^{i+1}}{k!} = C \frac{x^{n+2}}{(n+1)!}, \]
and 
\[ |\var{T} - n(1+W(n))| \leq \sum_{i=n+1}^\infty \frac{x^i}{i!}. \]

With respect to Theorem~\ref{p_qlclt_for_combinatorics}, specifically the set $M$, we also need to compute the maximum among the variances of $i\, Z_i$.
We define 
\[ b_i := \var{i\, Z_i} = i^2 \frac{x^i}{i!}. \]
One can appeal to the heuristics in~\cite{IPARCS} to analyze this part.  
The first few variances are
\[ b_1 = x \sim \log(n), \ \ b_2 = \frac{4x^2}{2} \sim 2\log^2(n), \ \ b_3 = \frac{9x^3}{6} \sim \frac{3}{2} \log^3(n), \ \ \cdots \ . \]
Thus, for a fixed, finite positive integer $k$, $b_1, \ldots, b_k$ is an increasing sequence for $n$ large. 
This is in contrast to integer partitions, where the variances were asymptotically equivalent for part sizes of size $o(\sqrt{n})$, and within a constant factor for part sizes $O(\sqrt{n})$. 
Once $k$ grows large enough relative to $n$, however, the denominator dominates; this occurs when $x^i$ and $i!$ are roughly of the same order. 
Stirling's formula suggests something on the order of $i=[x]$. 
More precisely, for any real $\alpha$, $i = [x + \alpha\sqrt{x}]$ is the range of values of $i$ which are most relevant. 
From this, we obtain 
\[ b_{[x+\alpha\sqrt{x}]} \sim \frac{x^2\, x^x x^{\alpha \sqrt{n}}}{x^{x+\alpha\sqrt{x}} e^{-x-\alpha\sqrt{x}}(1+\alpha/\sqrt{x})^{x+\alpha\sqrt{x}} \sqrt{2\pi x}} \sim e^x\, \sqrt{\frac{x^3}{2\pi}}\ e^{-\alpha^2} \sim \frac{n \sqrt{\log{n}}}{\sqrt{2\pi}}\ e^{-\alpha^2}. \]
And so what we see is that indeed $\alpha = 0$ maximizes this variance, and that we may take the set $M$ in Theorem~\ref{p_qlclt_for_combinatorics} to be $\{[x+\alpha_1\sqrt{x}], \ldots, [x+\alpha_2\sqrt{x}]\}$ for any fixed, real $\alpha_1< \alpha_2$. 

However, this means that $e^{-|M|} = O\left(e^{-\sqrt{\log(n)}}\right)$, which decays slower than any positive power of $n$. 
In particular, we have $\sigma \sim \sqrt{n \log{n}}$, and hence $\sigma e^{-|M|}$ is not tending to zero, even though we have $\sigma_{\max}/\sigma = O(\log^{-1/4}(n)).$

\section{Proofs of technical lemmas}
\label{section:proofs}

\subsection{Proof of auxiliary results}\label{s_proof_of_auxiliary_results}

In this section we prove the auxiliary statements that were used in the proof of Theorem~\ref{p_qlclt_for_combinatorics}.\\

We start with the proof of Lemma~\ref{p_estimating_higher_moments_by_lower_moments}. We will follow a general strategy for proving Reverse H\"older inequalities for real-valued random variables that goes back to C.~Borell~\cite{borell1973inverse} (see also Barvinok's lecture notes~\cite[Theorem 26.1]{barvinokmeasure}). 
In this article, we adapt the argument to integer-valued random variables. The argument for real-valued random variables consists of two steps. In the first step one uses the Brunn-Minkowski inequality to show that the tails of a log-concave probability measure decay exponentially. In the second step, one uses the decay of the tails to establish the reverse H\"older inequality.\\

In our argument, we follow this general outline. It seems that one could reprove every step in the argument for integer-valued variables. However, to make the argument more interesting we will give a new alternative proof for the first step not relying on the Brunn-Minkowski inequality. We formulate the first step in the following auxiliary lemma.
\begin{lemma}\label{p_expon_tails_log_concave}
   Assume that~$X \in \mathbb{N}$ is a log-concave random variable i.e. for any~$k \in \mathbb{N}$ it holds
   \begin{align}\label{e_def_log_concave}
     \mathbb{P} \left[  X= k+1 \right]      \mathbb{P} \left[  X= k-1 \right] =      \mathbb{P} \left[  X= k \right]^2.
   \end{align}
We assume that the first moment of~$X$ is bounded i.e.
\begin{align}\label{e_revers_hoelder_aux_lemma_uniform_moment_bound}
  \mathbb{E} \left[ |X|\right] \leq C_1.
\end{align}
Then there are constants~$0<k_0< \infty$ and $\lambda <1$, that only depends on the constant~$C_1$ in~\eqref{e_revers_hoelder_aux_lemma_uniform_moment_bound}, such that for all~$k_0\leq k \in \mathbb{N}$
\begin{align}\label{e_exp_tails_log_concave}
  \mathbb{P} \left[ X =k \right] \leq  \lambda^k.
\end{align}
\end{lemma}
We first comment on Lemma~\ref{p_expon_tails_log_concave}. The statement~\eqref{e_exp_tails_log_concave} without the control of~$\lambda$ is well known in the literature (see for example~\cite[Section~4]{saumard2014log}). However, in order to carry out the second step one needs a uniform bound on the decay rate~$\lambda$. We want to note that in Lemma~\ref{p_expon_tails_log_concave} we give quite a precise control on the constant~$\lambda$: The constant~$\lambda$ only depends on the constant~$C_1$ bounding the first moment. We establish this uniform control on~$\lambda$ by using a similar indirect argument as in the proof of Lemma~\ref{p_control_characteristic_fnction_intermediate_and_large_values}. 

\begin{proof}[Proof of Lemma~\ref{p_expon_tails_log_concave}]
We start with showing that for any log-concave random variable~$X \in \mathbb{N}$ there are constants~$0< k_X < \infty$ and~$0< \lambda < 1$ such that for all~$k_X \leq k \in \mathbb{N}$
\begin{align}\label{e_individual_tail_bound_log_concave_RV}
  \mathbb{P} \left[X= k \right] \leq \lambda_X^k. 
\end{align}
For deducing~\eqref{e_individual_tail_bound_log_concave_RV} we follow the standard argument of~\cite[Section~4]{saumard2014log}. It follows directly from the definition of log-concavity~\eqref{e_def_log_concave} that for all~$k \geq 2$
\begin{align}
  \frac{\mathbb{P}[X=k]}{\mathbb{P}[X=k-1]} \geq \frac{\mathbb{P}[X=k+1]}{\mathbb{P}[X=k]}.
\end{align}
This means that the sequence~$0\geq \frac{\mathbb{P}[X=k]}{\mathbb{P}[X=k-1]}$ is decreasing in~$k$ and therefore the limit 
\begin{align}
  \lim_{k \to \infty} \frac{\mathbb{P}[X=k]}{\mathbb{P}[X=k-1]} = \theta
\end{align}
exists. We will now show that~$0< \theta < 1$, which then yields the desired statement by writing
\begin{align}
  \mathbb{P}[X=k] = \frac{\mathbb{P}[X=k]}{\mathbb{P}[X=k-1]} \frac{\mathbb{P}[X=k-1]}{\mathbb{P}[X=k-2]} \ldots \mathbb{P}[X=l].
\end{align}
In the last identity we assumed without loss of generality that~$\mathbb{P}[X=\tilde k]\neq 0 $ for~$l \leq \tilde k \leq k$.\\
Indeed, because log-concave random variables are unimodal it follows that
\begin{align}
  \lim_{k \to \infty} \frac{\mathbb{P}[X=k]}{\mathbb{P}[X=k-1]} = \theta \leq 1. 
\end{align}
Additionally, because~$X$ is log-concave it follows that~$\mathbb{E} \left[ |X| \right] < \infty $. Hence, we get that~$\theta < 1$ because else~$\mathbb{E} \left[ |X| \right] = \infty $.

Now, let us show that the constants~$k_X$ and~$\lambda_X$ only depend on the upper bound in~\eqref{e_revers_hoelder_aux_lemma_uniform_moment_bound}. The argument is indirect, so let us assume that this is not the case. Hence, there is a sequence of log-concave random variables~$X_l \in \mathbb{N}$ and integers~$k_l$ satisfying the uniform bound~\eqref{e_revers_hoelder_aux_lemma_uniform_moment_bound} such that~$k_l \to \infty$ and
\begin{align}\label{e_indirect_assumption}
\lim_{l \to \infty}  \mathbb{P}[X_l = k_l] =1.
\end{align}
Because the first moment of the random variables~$X_l$ is uniformly bounded it follows that the sequence of random variabels~$X_l$ is tight. Hence, there exists a subsequence that converges weakly to a random variable~$X_\infty$. This means that we can assume WLOG~that 
\begin{align}
  X_l \overset{w}{\rightarrow} X_\infty.
\end{align}
It follows from~\eqref{e_indirect_assumption} and the definition of weak convergence that
\begin{align}
  \mathbb{E} \left[ |X_\infty| \right] = \lim_{l \to \infty} \mathbb{E} \left[ |X_l| \right] \geq \liminf_{l \to \infty} k_l \mathbb{P}[X_l = k_l] = \infty, 
\end{align}
which means that~$X_\infty$ is not integrable. However, because weak convergence preserves log-concavity we know that~$X_{\infty}$ is log-concave, which means in particular that its first moment exists. This a contradiction.
\end{proof}
With the help of Lemma~\ref{p_expon_tails_log_concave}, we can now prove Lemma~\ref{p_estimating_higher_moments_by_lower_moments} as in Barvinok's lecture notes~\cite[Theorem 26.1]{barvinokmeasure}. 
\begin{proof}[Proof of Lemma~\ref{p_estimating_higher_moments_by_lower_moments}.]
  It suffices to show that there is a universal constant such that
  \begin{align*}
    \left( \mathbb{E} \left[ |X|^p \right]\right)^{\frac{1}{p}} \leq C C_1^{\frac{1}{p}+ \frac{1}{q}} p \mathbb{E} \left[ |X| \right] .
  \end{align*}
Using that the random variable~$X$ is perurbed log-concave in the sense of~\eqref{e_def_pertubed_log_concave} we can assume WLOG~ that~$X$ is log-concave and it suffices to show that
  \begin{align}\label{e_reduced_estimate_reverse_hoelder}
    \left( \mathbb{E} \left[ |X|^p \right]\right)^{\frac{1}{p}} \leq C p \mathbb{E} \left[ |X| \right] .
  \end{align}
Additionally, we may assume WLOG~that
\begin{align}\label{e_proof_reverse_Hoelder_wlog}
  \mathbb{E} \left[ |X| \right]= 1.
\end{align}
Using the layer-cake presentation we get 
\begin{align} \label{e_layer_cake_representation}
  \mathbb{E} \left[ |X|^p \right] & = \int_0^\infty t^p d F(t),
\end{align}
where~$F(t)$ is the distribution fuction of~$X$ i.e.
\begin{align}
  F(t) = \mathbb{P} \left( X \leq t\right).
\end{align}
From Lemma~\ref{p_expon_tails_log_concave} we know that for~$k_0 \in \mathbb{N}$ large enough there is a unifrom constant~$0<c< 1$ such that
\begin{align}\label{e_tail_estimate_distribuition_function}
  F(t) = \mathbb{P} \left( X \leq t\right) \leq \exp (- c t).
\end{align}
Using the estimate~\eqref{e_tail_estimate_distribuition_function} we can estimate~\eqref{e_layer_cake_representation} as
\begin{align*}
  \mathbb{E} \left[ |X|^p \right] & \leq \int_0^{k_0} t^p dt + \int_{k_0}^\infty t^p d F(t)   \\
& \leq \frac{1}{p+1} k_0^{p+1} + \int_{k_0}^\infty t^p \exp(-ct) dt\\
& \leq \frac{1}{p+1} k_0^{p+1} + \frac{1}{c^{p+1}} \int_{0}^\infty t^p \exp(-t) dt\\
& \leq \frac{1}{p+1} k_0^{p+1} + \frac{1}{c^{p+1}} \Gamma(p),
\end{align*}
which yields the desired estimate~\eqref{e_reduced_estimate_reverse_hoelder} by using~\eqref{e_proof_reverse_Hoelder_wlog} and that the Gamma function statisfies~$\Gamma(p) \leq p^p$.
\end{proof}

\begin{proof}[Proof of Lemma~\ref{p_taylor_characteristic_function}]
The idea of the argument is to deduce the desired estimate~\eqref{e_control_characteristic_taylor} using Taylor expansion.
Therefore let us calculate the first three derivatives of the function~$h (\xi)$. The first derivative of~$h$ is given by
\begin{align}\label{e_first_derivative_h}
  \frac{d}{d\xi} h(\xi) = \frac{\mathbb{E} [i X \exp(i X \xi)]}{\mathbb{E} [\exp(i X \xi)]} ,
\end{align}
the second derivative is given by
\begin{align}
  \frac{d^2}{d\xi^2} h(\xi) = \frac{\mathbb{E} \left[ \left( iX - \frac{\mathbb{E}[iX\exp(i X \xi)]}{\mathbb{E}[\exp(i X \xi)]}\right)^2 \exp(i X \xi)\right]}{\mathbb{E} [\exp(i X \xi)]}, \label{e_second_derivative_h}
\end{align}
and the third derivative is given by
\begin{align}
  \frac{d^3}{d\xi^3} h(\xi) & = -2\ \frac{\mathbb{E} \left[ \left( iX - \frac{\mathbb{E}[iX\exp(i X \xi)]}{\mathbb{E}[\exp(i X \xi)]}\right) \frac{\mathbb{E} \left[ \left( iX - \frac{\mathbb{E}[iX\exp(i X \xi)]}{\mathbb{E}[\exp(i X \xi)]}\right)^2 \exp(i X \xi)\right]}{\mathbb{E} [\exp(i X \xi)]} \exp(i X \xi)\right]}{\mathbb{E} [\exp(i X \xi)]} \\
   & \qquad + \frac{\mathbb{E} \left[ \left( iX - \frac{\mathbb{E}[iX\exp(i X \xi)]}{\mathbb{E}[\exp(i X \xi)]}\right)^2 \left( iX - \frac{\mathbb{E}[iX\exp(i X \xi)]}{\mathbb{E}[\exp(i X \xi)]}\right)  \exp(i X \xi)\right]}{\mathbb{E} [\exp(i X \xi)]}. \label{e_third_derivative_h}
\end{align}

 Direct calculation yields that
  \begin{align*}
   \frac{d}{d\xi} h(0)= \mathbb{E}[X]=0 \qquad \mbox{and} \qquad \frac{d^2}{d\xi^2} h(0) = \var(X) =1.
  \end{align*}
Therefore in order to deduce~\eqref{e_control_characteristic_taylor} it suffices to show that for~$|\xi|\leq \delta$ one has a uniform bound
\begin{align*}
\left|\frac{d^3}{d\xi^3} h(\xi) \right| \leq C < \infty.
\end{align*}
Straightforward estimation using H\"older's inequality yields the estimate
  \begin{align*}
  \left| \frac{d^3}{d\xi^3} h(\xi) \right| & \leq C\, \mathbb{E} |X|^3 \left( \frac{1} {\left| \mathbb{E} [\exp(i X \xi)]\right|}  + \frac{1} {\left| \mathbb{E} [\exp(i X \xi)]\right|^2} + \frac{1} {\left| \mathbb{E} [\exp(i X \xi)]\right|^3} \right).
\end{align*}
The argument is finished if we show that by choosing~$\delta$ small enough one gets for~$|\xi| \leq \delta$
\begin{align*}
\left| \mathbb{E} [\exp(i X \xi)]\right| \geq \frac{1}{2}.
\end{align*}
We use again Taylor series expansion. Setting~$g(\xi) = \mathbb{E} [\exp(i X \xi)]$ we get by direct calculation that
\begin{align*}
\frac{d}{d \xi} g(\xi) =\frac{d}{d\xi}  \mathbb{E} [\exp(i X \xi)] = \mathbb{E} \left[ iX \exp(iX \xi) \right],
\end{align*}
which implies by our normalization~$\mathbb{E} [X] =0$ that
\begin{align*}
  \frac{d}{d\xi} g(0) =0.
\end{align*}
By the normalization~$\var(X) =1$ we have that
\begin{align*}
\left| \frac{d}{d\xi}  \mathbb{E} [\exp(i X \xi)] \right|  \leq \mathbb{E} \left|X\right| \leq 1.
\end{align*}
This together with the observation that~$g(0) = 1$ implies by Taylor expansion that for~$|\xi| \leq \delta$
\begin{align*}
  \left|  g(\xi) -1 \right| \leq  |\xi| \leq \delta .
\end{align*}
Therefore we obtain the estimate
\begin{align*}
  |g(\xi)| = | 1+ g(\xi) -1| \geq 1 -  |g(\xi) -1| \geq (1- \delta) \geq \frac{1}{2},
\end{align*}
which closes the argument.
\end{proof}

\begin{proof}[Proof of Lemma~\ref{p_control_characteristic_fnction_intermediate_and_large_values}]
  Let us assume the opposite statement of our Lemma.
Then there exists a sequence of random variables~$X_n$ and numbers~$\delta\leq |\xi_n| \leq 2\pi -\delta$ such that for all~$n\in \mathbb{N}$ the following five conditions are satisfied:

\begin{itemize}
 \item We have the uniform upper bound

\begin{align}\label{e_tightness_random_variables}
  \mathbb{E} \left[ |X_n|\right]\leq C_2< \infty.
\end{align}
\item We have the uniform lower bound
  \begin{align}\label{e_variance_random_variables}
      \var (X_n) \geq \frac{1}{C_3}> 0.
  \end{align}
\item The random variables~$X_n$ take values on a lattice i.e.
\begin{align}\label{e_lattice_random_variable}
  \mathbb{P} \left( X_n= c_{1,n} k + c_{2,n}, \ k \in \mathbb{Z}  \right) =1 \qquad \mbox{for constants } 0 \leq c_{1,n},c_{2,n} \leq 1
\end{align}
\item The probability distribution function of~$X_n$ is perturbed log-concave (see Definition~\ref{perturbed:unimodal}) in the sense that there are log-concave random variables~$Y_n$ such that
    \begin{align}\label{e_bounded_perturbation_unimodal}
\frac{1}{C_4}\leq \frac{\mathbb{P}(X_n=x)  }{\mathbb{P}(Y_n=x)} \leq C_4 \qquad \mbox{for all } x \in \mathbb{R},
\end{align}
where we define~$\frac{0}{0}=0$.

\item The characteristic functions evaluated at~$\xi_n$ satisfy
\begin{align*}
  \lim_{n \to \infty} |\mathbb{E} \left[ \exp (i X_n \xi_n) \right]|  \to 1.
\end{align*}
\end{itemize}
\mbox{}\medskip

By using that the numbers~$\xi_n$,~$c_{1,n}$ and~$c_{2,n}$ are uniformly bounded, we may choose a subsequence such that those numbers converge to some limit. More precisely, we may assume without loss of generality that there exists numbers~$\delta \leq |\xi_{\infty}|\leq \pi, 0 \leq c_{1,\infty}, c_{2,\infty}  \leq 1$ such that
\begin{align*}
\lim_{n \to \infty} \xi_{n} = \xi_{\infty},  \quad \lim_{n \to \infty} c_{1,n} = c_{1, \infty} \quad \lim_{n \to \infty} c_{2,n} = c_{2, \infty}.
\end{align*}
Because the random variables~$X_n$ satisfy the bound~\eqref{e_tightness_random_variables}, the family~$\left\{X_n \right\}$ is tight. Hence, it follows from Prokhorov's theorem that there exists a further subsequence and a random variable~$X_{\infty}$ such that~$X_n$ converges weakly to~$X_\infty$. This means that for all bounded measurable and continous functions~$f$ 
\begin{align*}
  \mathbb{E}\left[ f(X_n)\right] \to \mathbb{E}\left[ f(X_{\infty})\right] .
\end{align*}
This implies in particular that by~\eqref{e_variance_random_variables}
\begin{align}\label{e_variance_lower_bound_X_infty}
  \var(X_\infty) = \lim_{n \to \infty} \var( X_{n}) \geq \frac{1}{C_3}
\end{align}
We can choose also the function~$f_{\varepsilon}$ as
\begin{align*}
  f_{\varepsilon}(x) = \sum_{k \in \mathbb{Z}} \mathds{1}_{ \left\{ |c_{1,\infty} x - k + c_{2, \infty}| \leq \varepsilon \right\} } (x).
\end{align*}
By the weak convergence and~\eqref{e_lattice_random_variable} we get that for all~$\varepsilon> 0$
\begin{align*}
  \mathbb{E}\left[f_{\varepsilon}(X_{\infty}) \right] = \lim_{n \to \infty} \mathbb{E}\left[f_{\varepsilon}(X_{n}) \right] =1.
\end{align*}
Sending~$\varepsilon \to 0$, the last identity implies that
\begin{align}\label{e_values_X_infty}
  \mathbb{P} \left( X_\infty= c_{1,\infty} k + c_{2,\infty}, \ k \in \mathbb{Z}  \right) =1.
\end{align}
By the estimate~\eqref{e_variance_lower_bound_X_infty}, we know that the random variable~$X_\infty$ is not degenerate. This implies that
\begin{align}\label{e_X_infty_non_degenerate}
  c_{1,\infty} \neq 0.
\end{align}
By construction it also holds that
\begin{align*}
 \left|  \mathbb{E} \left[\exp \left(i X_\infty \xi_\infty \right) \right] \right|=1,
\end{align*}
hence there is~$\zeta \in S^1$ such that 
\begin{align}\label{e_pointwise_identity}
  \exp \left(i X_\infty \xi_\infty \right) = \zeta \qquad \mbox{for almost every} ~X_\infty.
\end{align}
Because of the bound~\eqref{e_variance_lower_bound_X_infty} we know that~$X_{\infty}$ takes on at least two values~$x_1$ and~$x_2$ with non-zero probabilities. By~\eqref{e_values_X_infty}, we can find~$k \neq \tilde k  \in \mathbb{Z}$ such that
\begin{align*}
  x_1= c_{1,\infty} k + c_{2,\infty} \qquad \mbox{and} \qquad x_2=c_{1,\infty} \tilde k + c_{2,\infty} .
\end{align*}
Hence, we get that
\begin{align}\label{e_nonegative_probability_at_k_and_tilde_k}
  \mathbb{P} (X_\infty=x_1) \neq 0  \quad \mbox{and} \quad  \mathbb{P} (X_\infty =x_2 ) \neq 0.
\end{align}
 \medskip

We didn't use yet the property that the random variables~$X_n$ are bounded perturbation of log-concave random variables. We will use this property to show that if for any~$k< \tilde k \in \mathbb{Z}$
\begin{align*}
  \mathbb{P} (X_\infty=x_1=c_{1,\infty} k + c_{2,\infty}) \neq 0  \quad \mbox{and} \quad  \mathbb{P} (X_\infty =x_2= c_{1,\infty} \tilde k + c_{2,\infty} ) \neq 0
\end{align*}
then it also holds that for~$k+1$
\begin{align} \label{e_nonegative_probability_at_k_plus_1}
    \mathbb{P} (X_\infty= c_{1,\infty} (k+1) + c_{2,\infty}) \neq 0 .
\end{align}
Let us assume for the moment, that this conclusion is true. Then we get by~\eqref{e_pointwise_identity},~\eqref{e_nonegative_probability_at_k_and_tilde_k} and~\eqref{e_nonegative_probability_at_k_plus_1} that 
\begin{align*}
  \exp \left(i \left( c_{1,\infty} k + c_{2,\infty} \right) \xi_\infty \right) & = \exp \left(i x_1 \xi_\infty \right)  \\
                                                                              & = \zeta \\
  & = \exp \left(i x_2 \xi_\infty \right) \\ 
  & = \exp \left(i \left(   c_{1,\infty} (k+1) + c_{2,\infty} \right) \xi_\infty \right).
\end{align*}
The last identity implies
\begin{align*}
  \exp \left(i  c_{1,\infty} \xi_\infty \right) = 1,
\end{align*}
which implies that~$ c_{1,\infty} \xi_\infty = m 2 \pi $ for some~$m \in \mathbb{Z}$. Because~$ c_{1,\infty}\neq 0$ by~\eqref{e_X_infty_non_degenerate} and~$|\xi_\infty|\geq \delta$ by construction, it follows that~$ |c_{1,\infty} \xi_\infty| \geq 2 \pi $. Using that~$|\xi_\infty| \leq 2\pi - \delta$ and $|c_{1,\infty}| \leq 1 $ it follows that~$| c_{1,\infty} \xi_\infty| \leq 2 \pi- \delta$ which is contradiction to the statement~$ |c_{1,\infty} \xi_\infty| \geq 2 \pi $ and therefore concludes the argument.\medskip

It is only left to show that the conclusion~\eqref{e_nonegative_probability_at_k_plus_1} is true. As mentioned before we will use the fact that the random variables~$X_n$ are perturbed log-concave. Because of~\eqref{e_tightness_random_variables} and~\eqref{e_bounded_perturbation_unimodal}, we know that the first moment of the log-concave random variables~$Y_n$ are uniformly bounded. More precisely, we can show that
\begin{align}\label{e_tightness_uninomdal_Y_n}
  \mathbb{E} \left[ |Y_n| \right] \leq C_3 C_4.
\end{align}
It follows from~\eqref{e_tightness_uninomdal_Y_n} that the family of random variables~$Y_n$ is also tight. Hence there exists a subsequence such that~$Y_n $ converges weakly to a random variable~$Y_\infty$. With the same reasoning as for the random variable~$X_\infty$, a combination of~\eqref{e_lattice_random_variable} and~\eqref{e_bounded_perturbation_unimodal} yields that the random variable~$Y_\infty$ is also a lattice random variable with parameters~$c_{1,\infty}$ and~$c_{2 ,\infty}$; i.e.,
\begin{align}\label{e_Y_infty_lattice}
    \mathbb{P} \left( Y_\infty= c_{1,\infty} k + c_{2,\infty}, \ k \in \mathbb{Z}  \right) =1.
\end{align}
Because~$X_\infty$ and $Y_\infty$  are lattice random variables, it follows from the weak convergence of~$X_n \Rightarrow X_\infty$ and~$Y_n \Rightarrow Y_\infty$ that for all~$x\in \mathbb{R}$
\begin{align}\label{e_convergence_probabilities_of_X_n_and_Y_n}
 \lim_{n \to \infty} \mathbb{P} (X_n =x) = \mathbb{P}( X_\infty=x) \quad \mbox{and}  \quad \lim_{n \to \infty} \mathbb{P} (Y_n =x)= \mathbb{P}(Y_\infty=x).
\end{align}
Using now the property~\eqref{e_Y_infty_lattice} allows us to conclude that the random variable~$Y_\infty$ inherits the log-concavity from the random variables~$Y_n$. A combination of~\eqref{e_bounded_perturbation_unimodal} and~\eqref{e_convergence_probabilities_of_X_n_and_Y_n} yields that
\begin{align}\label{e_bounded_perturbation_unimodal_limit}
 \frac{1}{C_4}\leq \frac{\mathbb{P}(X_\infty=x)  }{\mathbb{P}(Y_\infty=x)} \leq C_4 \qquad \mbox{for all } x \in \mathbb{R} .
\end{align}
However, because~$Y_\infty$ is log-concave, the last identity yields that the probability distribution function of the random variable~$X_\infty$ is a bounded perturbation of the a log-concave probability distribution function of the random variable~$Y_\infty$. 
Now, let us assume that for~$k < \tilde k \in \mathbb{Z}$ it holds
\begin{align*}
  \mathbb{P} (X_\infty=x_1=c_{1,\infty} k + c_{2,\infty}) \neq 0  \quad \mbox{and} \quad  \mathbb{P} (X_\infty =x_2= c_{1,\infty} \tilde k + c_{2,\infty} ) \neq 0.
\end{align*}
By the estimate~\eqref{e_bounded_perturbation_unimodal_limit} it follows that
\begin{align}\label{e_nonnegative_probabilities_Y_infty}
  \mathbb{P} (Y_\infty=x_1=c_{1,\infty} k + c_{2,\infty}) \neq 0  \quad \mbox{and} \quad  \mathbb{P} (Y_\infty =x_2= c_{1,\infty} \tilde k + c_{2,\infty} ) \neq 0.
\end{align}
A combination of~\eqref{e_nonnegative_probabilities_Y_infty} and the fact that the random variable~$Y_\infty$ is unimodal yields that for all~$\hat k \in \mathbb{Z}$,~$k \leq \hat k  \leq \tilde k$
\begin{align*}
   \mathbb{P} (Y_\infty=c_{1,\infty} \hat k + c_{2,\infty}) \neq 0.
\end{align*}
In particular by setting~$\hat k = k+1$ this yields that
\begin{align*}
   \mathbb{P} (Y_\infty=c_{1,\infty} (k+1) + c_{2,\infty}) \neq 0.
\end{align*}
Using now the estimate~\eqref{e_bounded_perturbation_unimodal_limit} finally yields the desired statement that
\begin{align*}
   \mathbb{P} (X_\infty=c_{1,\infty} (k+1) + c_{2,\infty}) \neq 0. 
\end{align*}
\end{proof}

\bibliographystyle{alpha}
\bibliography{bib_one}

\newcommand{\etalchar}[1]{$^{#1}$}
\begin{thebibliography}{CCH01b}

\bibitem[ABT03]{Logarithmic}
Richard Arratia, Andrew~D Barbour, and Simon Tavar{\'e}.
\newblock {\em Logarithmic combinatorial structures: a probabilistic approach},
  volume~1.
\newblock European Mathematical Society Z{\"u}rich, 2003.

\bibitem[AD16]{PDC}
Richard Arratia and Stephen DeSalvo.
\newblock Probabilistic divide-and-conquer: a new exact simulation method, with
  integer partitions as an example.
\newblock {\em Combinatorics, Probability and Computing}, 25(3):324--351, May
  2016.

\bibitem[Alm02]{almkvist2002partitions}
Gert Almkvist.
\newblock Partitions with parts in a finite set and with parts outside a finite
  set.
\newblock {\em Experiment. Math.}, 11(4):449--456 (2003), 2002.

\bibitem[And84]{Andrews}
George~E. Andrews.
\newblock {\em The Theory of Partitions}.
\newblock Cambridge Mathematical Library, 1984.

\bibitem[And00]{AndrewsII}
George~E. Andrews.
\newblock Mac{M}ahon's partition analysis. {II}. {F}undamental theorems.
\newblock {\em Ann. Comb.}, 4(3-4):327--338, 2000.
\newblock Conference on Combinatorics and Physics (Los Alamos, NM, 1998).

\bibitem[APR01]{AndrewsOmega}
George~E. Andrews, Peter Paule, and Axel Riese.
\newblock Mac{M}ahon's partition analysis: the {O}mega package.
\newblock {\em European J. Combin.}, 22(7):887--904, 2001.

\bibitem[AST95]{arratia1995total}
Richard Arratia, Dudley Stark, and Simon Tavar{{\'e}}.
\newblock Total variation asymptotics for {P}oisson process approximations of
  logarithmic combinatorial assemblies.
\newblock {\em Ann. Probab.}, 23(3):1347--1388, 1995.

\bibitem[AT94]{IPARCS}
Richard Arratia and Simon Tavare.
\newblock Independent process approximations for random combinatorial
  structures.
\newblock {\em Adv. Math. 104 (1994), no. 1, 90-154}, 08 1994.

\bibitem[Bar]{barvinokmeasure}
A~Barvinok.
\newblock Measure concentration lecture notes, 2005.

\bibitem[B{\'C}{\etalchar{+}}02]{Barbour}
Andrew~D Barbour, V~{\'C}ekanavi{\'c}ius, et~al.
\newblock Total variation asymptotics for sums of independent integer random
  variables.
\newblock {\em The Annals of Probability}, 30(2):509--545, 2002.

\bibitem[BE56]{ErdosBateman}
P.~T. Bateman and P.~Erd{\"o}s.
\newblock Monotonicity of partition functions.
\newblock {\em Mathematika}, 3:1--14, 1956.

\bibitem[Ben73]{Bender}
Edward~A. Bender.
\newblock Central and local limit theorems applied to asymptotic enumeration.
\newblock {\em J. Combinatorial Theory Ser. A}, 15:91--111, 1973.

\bibitem[BFP10]{BodiniPlane}
Olivier Bodini, {\'E}ric Fusy, and Carine Pivoteau.
\newblock Random sampling of plane partitions.
\newblock {\em Combinatorics, Probability and Computing}, 19(02):201--226,
  2010.

\bibitem[BM11]{BoMA11}
Sergey Bobkov and Mokshay Madiman.
\newblock Concentration of the information in data with log-concave
  distributions.
\newblock {\em Ann. Probab.}, 39(4):1528--1543, 2011.

\bibitem[Bog15]{Bogachev}
Leonid~V. Bogachev.
\newblock Unified derivation of the limit shape for multiplicative ensembles of
  random integer partitions with equiweighted parts.
\newblock {\em Random Structures Algorithms}, 47(2):227--266, 2015.

\bibitem[Bor73]{borell1973inverse}
Christer Borell.
\newblock Inverse h{\"o}lder inequalities in one and several dimensions.
\newblock {\em Journal of Mathematical Analysis and Applications},
  41(2):300--312, 1973.

\bibitem[CCH01a]{NonNegR}
Rod Canfield, Sylvie Corteel, and Pawel Hitczenko.
\newblock Random partitions with non-negative rth differences.
\newblock {\em Advances in Applied Mathematics}, 27(2):298--317, 2001.

\bibitem[CCH01b]{canfield2001random}
Rod Canfield, Sylvie Corteel, and Pawel Hitczenko.
\newblock Random partitions with non-negative rth differences.
\newblock {\em Advances in Applied Mathematics}, 27(2):298--317, 2001.

\bibitem[CGS10]{chen2010normal}
Louis~HY Chen, Larry Goldstein, and Qi-Man Shao.
\newblock {\em Normal approximation by Stein's method}.
\newblock Springer Science \&amp; Business Media, 2010.

\bibitem[Che75]{chen1975poisson}
Louis~HY Chen.
\newblock Poisson approximation for dependent trials.
\newblock {\em The Annals of Probability}, pages 534--545, 1975.

\bibitem[CW12]{CanfieldWilf}
E.~Rodney Canfield and Herbert~S. Wilf.
\newblock On the growth of restricted integer partition functions.
\newblock In {\em Partitions, {$q$}-series, and modular forms}, volume~23 of
  {\em Dev. Math.}, pages 39--46. Springer, New York, 2012.

\bibitem[DB70]{Debruijn}
Nicolaas~Govert De~Bruijn.
\newblock {\em Asymptotic methods in analysis}, volume~4.
\newblock Courier Dover Publications, 1970.

\bibitem[DLP16]{Dalal}
Avinash~J Dalal, Amanda Lohss, and Daniel Parry.
\newblock Statistical structures of concave compositions.
\newblock {\em arXiv preprint arXiv:1605.00343}, 2016.

\bibitem[DM95]{Davis}
Burgess Davis and David McDonald.
\newblock An elementary proof of the local central limit theorem.
\newblock {\em Journal of Theoretical Probability}, 8(3):693--701, 1995.

\bibitem[DN90a]{dixmier1990partitions}
J.~Dixmier and J.-L. Nicolas.
\newblock Partitions without small parts.
\newblock In {\em Number theory, {V}ol.\ {I} ({B}udapest, 1987)}, volume~51 of
  {\em Colloq. Math. Soc. J{\'a}nos Bolyai}, pages 9--33. North-Holland,
  Amsterdam, 1990.

\bibitem[DN90b]{dixmier1990partitions2}
Jacques Dixmier and Jean-Louis Nicolas.
\newblock Partitions sans petits sommants.
\newblock In {\em A tribute to {P}aul {E}rd{\H o}s}, pages 121--152. Cambridge
  Univ. Press, Cambridge, 1990.

\bibitem[DP13]{logconcave}
Stephen DeSalvo and Igor Pak.
\newblock Log-concavity of the partition function.
\newblock {\em The Ramanujan Journal}, 10 2013.

\bibitem[DVZ98]{LD}
Amir Dembo, Anatoly Vershik, and Ofer Zeitouni.
\newblock Large deviations for integer partitions, 1998.

\bibitem[EG08]{erlihson2008limit}
Michael~M. Erlihson and Boris~L. Granovsky.
\newblock Limit shapes of {G}ibbs distributions on the set of integer
  partitions: the expansive case.
\newblock {\em Ann. Inst. Henri Poincar{\'e} Probab. Stat.}, 44(5):915--945,
  2008.

\bibitem[EL41]{ErdosLehner}
Paul Erd{\H{o}}s and Joseph Lehner.
\newblock The distribution of the number of summands in the partitions of a
  positive integer.
\newblock {\em Duke Math. J}, 8(2):335--345, 1941.

\bibitem[ER76]{erdos1976concerning}
P.~Erd{\H{o}}s and B.~Richmond.
\newblock Concerning periodicity in the asymptotic behaviour of partition
  functions.
\newblock {\em J. Austral. Math. Soc. Ser. A}, 21(4):447--456, 1976.

\bibitem[Erd42]{ErdosElementary}
P.~Erd{{\"o}}s.
\newblock On an elementary proof of some asymptotic formulas in the theory of
  partitions.
\newblock {\em Ann. of Math. (2)}, 43:437--450, 1942.

\bibitem[Fri93]{Fristedt}
Bert Fristedt.
\newblock The structure of random partitions of large integers.
\newblock {\em Transactions of the American Mathematical Society},
  337(2):703--735, 1993.

\bibitem[FS09]{Flajolet}
Philippe Flajolet and Robert Sedgewick.
\newblock {\em Analytic combinatorics}.
\newblock Cambridge University Press, Cambridge, 2009.

\bibitem[GH08]{Goh}
William~MY Goh and Pawel Hitczenko.
\newblock Random partitions with restricted part sizes.
\newblock {\em Random Structures \&amp; Algorithms}, 32(4):440--462, 2008.

\bibitem[GK06]{grabner2006analysis}
Peter~J. Grabner and Arnold Knopfmacher.
\newblock Analysis of some new partition statistics.
\newblock {\em Ramanujan J.}, 12(3):439--454, 2006.

\bibitem[GS92]{goh1992gap}
William M.~Y. Goh and Eric Schmutz.
\newblock Gap-free set partitions.
\newblock {\em Random Structures Algorithms}, 3(1):9--18, 1992.

\bibitem[GS02]{LCLTGamkrelidzeLattice}
N.~Gamkrelidze and T.~Shervashidze.
\newblock On a local limit theorem for lattice distributions.
\newblock {\em Proc. A. Razmadze Math. Inst.}, 130:7--12, 2002.

\bibitem[Hag71]{hagis1971partitions}
Peter Hagis, Jr.
\newblock Partitions with a restriction on the multiplicity of the summands.
\newblock {\em Trans. Amer. Math. Soc.}, 155:375--384, 1971.

\bibitem[Hei90]{heinrich1990asymptotic}
Lothar Heinrich.
\newblock Asymptotic expansions in the central limit theorem for a special
  class of {$m$}-dependent random fields. {II}. {L}attice case.
\newblock {\em Math. Nachr.}, 145:309--327, 1990.

\bibitem[HR18]{HR}
G.H. Hardy and S.~Ramanujan.
\newblock Asymptotic formula{\ae} in combinatory analysis.
\newblock {\em Proceedings of the London Mathematical Society}, 2(1):75--115,
  1918.

\bibitem[Hua42]{hua1942number}
Loo-keng Hua.
\newblock On the number of partitions of a number into unequal parts.
\newblock {\em Trans. Amer. Math. Soc.}, 51:194--201, 1942.

\bibitem[Ing41]{Ingham}
AE~Ingham.
\newblock A tauberian theorem for partitions.
\newblock {\em Annals of Mathematics}, pages 1075--1090, 1941.

\bibitem[KR13]{kane2013asymptotics}
Daniel Kane and ROBERT~C Rhoades.
\newblock Asymptotics for wilf's partitions with distinct multiplicities.
\newblock {\em Preprint}, 2013.

\bibitem[Mac60]{MacMahon}
Percy~A. MacMahon.
\newblock {\em Combinatory analysis}.
\newblock Two volumes (bound as one). Chelsea Publishing Co., New York, 1960.

\bibitem[Mah40]{mahler1940special}
Kurt Mahler.
\newblock On a special functional equation.
\newblock {\em J. London Math. Soc.}, 15:115--123, 1940.

\bibitem[McD80]{McDonald}
David~R McDonald.
\newblock On local limit theorem for integer-valued random variables.
\newblock {\em Theory of Probability \&amp; Its Applications}, 24(3):613--619,
  1980.

\bibitem[Mei53]{Meinardus}
G{\"u}nter Meinardus.
\newblock Asymptotische aussagen {\"u}ber partitionen.
\newblock {\em Mathematische Zeitschrift}, 59(1):388--398, 1953.

\bibitem[Min72]{LCLTMineka}
J.~Mineka.
\newblock Local limit theorems and recurrence conditions for sums of
  independent integer-valued random variables.
\newblock {\em Ann. Math. Statist.}, 43:251--259, 1972.

\bibitem[MO13]{Menz_Otto}
Georg Menz and Felix Otto.
\newblock {U}niform logarithmic {S}obolev inequalities for conservative spin
  systems with super-quadratic single-site potential.
\newblock {\em Ann. Probab.}, 41(3B):2182--2224, 05 2013.

\bibitem[Muk91]{LCLTMukhin}
A.~B. Mukhin.
\newblock Local limit theorems for lattice random variables.
\newblock {\em Teor. Veroyatnost. i Primenen.}, 36(4):660--674, 1991.

\bibitem[Nat00]{nathanson2000partitions}
Melvyn~B. Nathanson.
\newblock Partitions with parts in a finite set.
\newblock {\em Proc. Amer. Math. Soc.}, 128(5):1269--1273, 2000.

\bibitem[New06]{newman2006analytic}
Donald~J Newman.
\newblock {\em Analytic number theory}, volume 177.
\newblock Springer Science \&amp; Business Media, 2006.

\bibitem[NS00]{nicolas2000partitions}
J.-L. Nicolas and A.~S{{\'a}}rk{{\"o}}zy.
\newblock On partitions without small parts.
\newblock {\em J. Th{\'e}or. Nombres Bordeaux}, 12(1):227--254, 2000.

\bibitem[Odl95]{odlyzko1995asymptotic}
A.~M. Odlyzko.
\newblock Asymptotic enumeration methods.
\newblock In {\em Handbook of combinatorics, {V}ol.\ 1,\ 2}, pages 1063--1229.
  Elsevier, Amsterdam, 1995.

\bibitem[Pak02]{PakBijection}
Igor Pak.
\newblock Hook length formula and geometric combinatorics.
\newblock {\em S{\'e}m. Lothar. Combin.}, 46:Art. B46f, 13 pp. (electronic),
  2001/02.

\bibitem[Pak04]{Pak2}
Igor Pak.
\newblock The nature of partition bijections ii.
\newblock {\em Asymptotic stability, preprint available at: http://www-math.
  mit. edu/pak/research. html}, 2004.

\bibitem[Pak06]{PakSurvey}
Igor Pak.
\newblock Partition bijections, a survey.
\newblock {\em The Ramanujan Journal}, 12(1):5--75, 2006.

\bibitem[Pet64]{LCLTPetrov}
Valentin~Vladimirovich Petrov.
\newblock On local limit theorems for sums of independent random variables.
\newblock {\em Theory of Probability \&amp; Its Applications}, 9(2):312--320,
  1964.

\bibitem[Pit97a]{PittelShape}
Boris Pittel.
\newblock On a likely shape of the random ferrers diagram.
\newblock {\em Adv. Appl. Math.}, 18(4):432--488, 1997.

\bibitem[Pit97b]{PittelSetPartitions}
Boris Pittel.
\newblock Random set partitions: asymptotics of subset counts.
\newblock {\em journal of combinatorial theory, Series A}, 79(2):326--359,
  1997.

\bibitem[Pit99]{PittelConfirming}
Boris Pittel.
\newblock Confirming two conjectures about the integer partitions.
\newblock {\em J. Combin. Theory Ser. A}, 88(1):123--135, 1999.

\bibitem[PW13]{pemantle2013analytic}
Robin Pemantle and Mark~C Wilson.
\newblock {\em Analytic combinatorics in several variables}, volume 140.
\newblock Cambridge University Press, 2013.

\bibitem[Ric75]{Richmond}
L.~B. Richmond.
\newblock Asymptotic relations for partitions.
\newblock {\em J. Number Theory}, 7(4):389--405, 1975.

\bibitem[Rob76]{robertson1976partitions}
M.~M. Robertson.
\newblock Partitions with congruence conditions.
\newblock {\em Proc. Amer. Math. Soc.}, 57(1):45--49, 1976.

\bibitem[R{\"o}l08]{Rollin}
Adrian R{\"o}llin.
\newblock Symmetric and centered binomial approximation of sums of locally
  dependent random variables.
\newblock {\em Electron. J. Prob}, 13:756--776, 2008.

\bibitem[Rom]{RomikUnpublished}
Dan Romik.
\newblock Identities arising from limit shapes of constrained random
  partitions.

\bibitem[Rom05]{Romik}
Dan Romik.
\newblock Partitions of $n$ into $t\sqrt{n}$ parts.
\newblock {\em Eur. J. Comb.}, 26(1):1--17, 2005.

\bibitem[Rom15]{romik2015number}
Dan Romik.
\newblock On the number of n-dimensional representations of su (3), the
  bernoulli numbers, and the witten zeta function.
\newblock {\em arXiv preprint arXiv:1503.03776}, 2015.

\bibitem[Roz57]{LCLTRozanov}
Yu.~A. Rozanov.
\newblock On a local limit theorem for lattice distributions.
\newblock {\em Teor. Veroyatnost. i Primenen.}, 2:275--281, 1957.

\bibitem[RS54]{RothSzekeres}
K.~F. Roth and G.~Szekeres.
\newblock Some asymptotic formulae in the theory of partitions.
\newblock {\em Quart. J. Math., Oxford Ser. (2)}, 5:241--259, 1954.

\bibitem[Sac74]{sachkov1974random}
Vladimir~Nikolaevich Sachkov.
\newblock Random partitions of sets.
\newblock {\em Theory of Probability \&amp; Its Applications}, 19(1):184--190,
  1974.

\bibitem[Ste72]{Stein}
Charles Stein.
\newblock A bound for the error in the normal approximation to the distribution
  of a sum of dependent random variables.
\newblock In {\em Proceedings of the {S}ixth {B}erkeley {S}ymposium on
  {M}athematical {S}tatistics and {P}robability ({U}niv. {C}alifornia,
  {B}erkeley, {C}alif., 1970/1971), {V}ol. {II}: {P}robability theory}, pages
  583--602. Univ. California Press, Berkeley, Calif., 1972.

\bibitem[SW14]{saumard2014log}
Adrien Saumard and Jon~A Wellner.
\newblock Log-concavity and strong log-concavity: a review.
\newblock {\em Statistics surveys}, 8:45, 2014.

\bibitem[Sze51]{Szekeres1}
George Szekeres.
\newblock An asymptotic formula in the theory of partitions.
\newblock {\em Quart. J. Math., Oxford Ser. (2)}, 2:85--108, 1951.

\bibitem[Sze53]{Szekeres2}
George Szekeres.
\newblock Some asymptotic formulae in the theory of partitions. {II}.
\newblock {\em Quart. J. Math., Oxford Ser. (2)}, 4:96--111, 1953.

\bibitem[Tem52]{Temperley}
Harold N.~V. Temperley.
\newblock Statistical mechanics and the partition of numbers ii. the form of
  crystal surfaces.
\newblock In {\em Mathematical Proceedings of the Cambridge Philosophical
  Society}, volume~48, pages 683--697. Cambridge Univ Press, 1952.

\bibitem[Ver96]{Vershik}
A.M. Vershik.
\newblock Statistical mechanics of combinatorial partitions, and their limit
  shapes.
\newblock {\em Functional Analysis and Its Applications}, 30:90--105, 1996.

\bibitem[VFY99]{Freiman}
Anatolii~Moiseevich Vershik, Gregory~A Freiman, and Yurii~Vladimirovich
  Yakubovich.
\newblock A local limit theorem for random strict partitions.
\newblock {\em Teoriya Veroyatnostei i ee Primeneniya}, 44(3):506--525, 1999.

\bibitem[VK77]{VershikKerov1977}
A.M. Vershik and S.V. Kerov.
\newblock Asymptotics of the {P}lancherel measure of the symmetric group and
  the limiting form of {Y}oung tableaux.
\newblock {\em Doklady Akademii Nauk SSSR}, 233(6):1024--1027, 1977.

\bibitem[Wri31]{WrightPlane}
E.~Maitland Wright.
\newblock Asymptotic partition formulae: {I}. plane partitions.
\newblock {\em Quarterly Journal of Mathematics, Oxford Series II}, pages
  177--189, 1931.

\bibitem[Wri34]{Wright}
E.~Maitland Wright.
\newblock Asymptotic partition formulae. {III}. {P}artitions into {$k$}-th
  powers.
\newblock {\em Acta Math.}, 63(1):143--191, 1934.

\bibitem[Yak12]{Yakubovich}
Yuri Yakubovich.
\newblock Ergodicity of multiplicative statistics.
\newblock {\em J. Comb. Theory Ser. A}, 119(6):1250--1279, August 2012.

\end{thebibliography}

\end{document}